\documentstyle[11pt]{article}
\def\arcsin{\mathop{\rm arcsin}\nolimits}

\def\arctan{\mathop{\rm arctan}\nolimits}

\def\erf{\mathop{\rm erf}\nolimits}

\newcommand{\C}{{\rm {\mbox{C{\llap{{\vrule height1.52ex}\kern.4em}}}}}}
\newcommand{\N} {{\rm {\mbox{\protect\makebox[.15em][l]{I}N}}}}
\newcommand{\Q} {{\rm {\mbox{Q{\llap{{\vrule height1.5ex}\kern.5em}}}}}}
\newcommand{\R} {{\rm {\mbox{\protect\makebox[.15em][l]{I}R}}}}

\newcommand{\Z} {{\rm {\mbox{\protect\makebox[.2em][l]{\sf Z}\sf Z}}}}

\newcommand{\ed}[1]{\frac{1}{#1}}

\newcommand{\funkdef}[3]{\left\{\begin{array}{ccc}
                                #1 && \mbox{\rm{if} $#2$ } \\
                                #3 && \mbox{\rm{otherwise}}
                                \end{array}  
                         \right.}       
\newcommand{\funkdeff}[4]{\left\{\begin{array}{ccc}
                                 #1 && \mbox{\rm{if} $#2$ } \\[3mm]
                                 #3 && \mbox{\rm{if} $#4$ } 
                                 \end{array}
                          \right.}
\newcommand{\funkdefff}[6]{\left\{\begin{array}{ccc}
                                 #1 && \mbox{\rm{if} $#2$ } \\
                                 #3 && \mbox{\rm{if} $#4$ } \\
                                 #5 && \mbox{\rm{if} $#6$ }
                                 \end{array}
                           \right.}
\newcommand{\funkdeffff}[8]{\left\{\begin{array}{ccc}
                                 #1 && \mbox{\rm{if} $#2$ } \\
                                 #3 && \mbox{\rm{if} $#4$ } \\
                                 #5 && \mbox{\rm{if} $#6$ } \\
                                 #7 && \mbox{\rm{if} $#8$ }
                                 \end{array}
                          \right.}
\newcommand{\ueber}[2]{
                       \Big( \!
                       {{\small
                       \begin{array}{c}
                          #1\\
                          #2
                          \end{array}
                       }}
                       \! \Big) }

\newcommand{\pr}{\vspace{-2mm}\absatz{{\sl Proof:}}\hspace{5mm}}
\newcommand{\eop}{\hfill$\Box$\par\medskip\noindent}
\newcommand{\absatz}{\par\medskip\noindent}

\newcommand{\al}{\alpha}

\newcommand{\la}{\lambda}

\renewcommand{\phi}{\varphi}

\newcommand{\th}{\theta}

\newcommand{\de}{\partial}

\newcommand{\sumi}{\sum\limits_{k=0}^{\infty}}

\newcommand{\1}{{\bf{1}}}
\newcommand{\2}{{\bf{2}}}
\newcommand{\3}{{\bf{3}}}

\newcommand{\5}{{\bf{5}}}
\newcommand{\6}{{\bf{6}}}
\newcommand{\7}{{\bf{7}}}

\newcommand{\pf}{\rightarrow}

\newcommand{\be}{\begin{equation}}
\newcommand{\ee}{\end{equation}}
\newcommand{\bea}{\begin{eqnarray}}
\newcommand{\eea}{\end{eqnarray}}
\newcommand{\beao}{\begin{eqnarray*}}
\newcommand{\eeao}{\end{eqnarray*}}

\newcommand{\lk}{\left(}
\newcommand{\rk}{\right)}

\newcommand{\rb}{\right|}
\newcommand{\bi}{\bibitem}

\newcommand{\bT}{\begin{theorem}}
\newcommand{\eT}{\end{theorem}}
\newcommand{\bL}{\begin{lemma}}
\newcommand{\eL}{\end{lemma}}
\newcommand{\bC}{\begin{corollary}}  
\newcommand{\eC}{\end{corollary}}
\newcommand{\bt}{\begin{tabbing} 12345 \= \kill}
\newcommand{\et}{\end{tabbing}}

\newtheorem{theorem}{Theorem}
\newtheorem{algorithm}{Algorithm}
\newtheorem{lemma}{Lemma}
\newtheorem{corollary}{Corollary}

\newcommand{\ded}[1]{\frac{\de}{\de #1}}
\newcommand{\dedn}[2]{\frac{\de^{#2}}{\de #1^{#2}}}
\def\airyai{\mathop{\rm Ai}\nolimits}
\def\airybi{\mathop{\rm Bi}\nolimits}
\def\erfc{\mathop{\rm erfc}\nolimits}
\def\erf{\mathop{\rm erf}\nolimits}

\begin{document}
\begin{center}{
{\LARGE {\bf {Algorithmic Work with}}}
\vspace{3mm}\\
{\LARGE {\bf {Orthogonal Polynomials}}}
\vspace{3mm}\\
{\LARGE {\bf {and Special Functions}}}
\vspace{5mm}\\
{\large {\sc {Wolfram Koepf}}}
\vspace{3mm}\\
{\sl Konrad-Zuse-Zentrum f\"ur Informationstechnik Berlin, Heilbronner Str. 10,
D-10711 Berlin, Federal Republic of Germany}
\\[3mm]
Konrad-Zuse-Zentrum Berlin (ZIB), Preprint SC 94-5, 1994
}
\end{center}
\vspace{0.5cm}
\begin{center}
{\bf {Abstract:}}
\end{center}
\begin{enumerate}
\item[]
{{\small
In this article we present a method to implement orthogonal polynomials
and many other special functions in Computer Algebra systems enabling
the user to work with those functions appropriately, and in particular
to verify different types of identities for those functions. Some of
these identities like differential equations, power series representations, 
and hypergeometric representations
can even dealt with algorithmically, i.\ e.\ they can
be computed by the Computer Algebra system, rather than only verified.

The types of functions that can be treated by the given technique
cover the generalized hypergeometric functions, and therefore most of the
special functions that can be found in mathematical dictionaries.

The types of identities for which we present verification algorithms cover
differential equations, power series representations, identities of the
Rodrigues type, hypergeometric representations, and algorithms containing
symbolic sums. 

The current implementations of special functions in existing Computer Algebra 
systems do not meet these high standards as we shall show in examples.
They should be modified, and we show results of our implementations.
}}
\end{enumerate}
%

\section{Introduction}
\label{sec:Introduction}

Many special functions
can be looked at from the following point of view: They represent
functions $f(n,x)$ of one ``discrete'' variable $n\in D$ defined on a set 
$D$ that has the property that $n\in D\Rightarrow n+1\in D$
(or $n\in D\Rightarrow n-1\in D$), e.\ g.\ $D=\N_0, \Z, \R$, or $\C$,
and one ``continuous'' variable $x\in I$ where $I$
represents a real interval, either finite $I=[a,b]$, infinite
($I=[a,\infty)$, $I=(-\infty,a]$, or $I=\R$), or a subset 
of the complex plane $\C$. 

In the given situation we may speak of the family 
$(f_n)_{n\in D}$ of functions $f_n(x):=f(n,x)$. 

In this paper we will deal with special functions and orthogonal polynomials 
of a real/complex variable $x$. Many of our results can be generalized to 
special and orthogonal functions of a discrete variable $x$ which we will
consider in a forthcoming paper.

Many of those families, especially all families of orthogonal polynomials, 
have the following properties:
\begin{enumerate}
\item
{\bf (Derivative rule)}
\\
The functions $f_n$ are differentiable with respect to the variable $x$,
and satisfy a derivative rule of the form
\be
f_n'(x)=\ded x f_n(x)=\sum_{k=0}^{m-1} r_k(n,x)\,f_{n-k}(x)
\quad\quad
\mbox{or}
\quad\quad
f_n'(x)=\sum_{k=0}^{m-1} r_k(n,x)\,f_{n+k}(x)
\;,
\label{eq:Derivative rule}
\ee
where the derivative with respect to $x$ is represented by a
finite number of lower or higher indexed functions of the family,
and where $r_k$ are rational functions in 
$x$. If $r_{m-1}(n,x)\not\equiv 0$ then the number 
$m$ is called the order of the given derivative rule.
We call the two different types of derivative rules {\sl backward} and
{\sl forward} derivative rule, respectively.
\item
{\bf (Differential equation)}
\\
The functions $f_n$ are $m$ times differentiable ($m\in\N)$ with respect
to the variable $x$, and satisfy a homogeneous
linear differential equation 
\be
\sum_{k=0}^m p_k(n,x)\,f_n^{(k)}(x)=0
\;,
\label{eq:Differential equation}
\ee
where $p_k$ are polynomials in 
$x$.
If $p_m(n,k)\not\equiv 0$ then the number $m$ is called the order
of the given differential equation.
\item
{\bf (Recurrence equation)}
\\
The functions $f_n$ satisfy a homogeneous
linear recurrence equation with respect to $n$
\be
\sum_{k=0}^m q_k(n,x)\,f_{n-k}(x)=0
\;,
\label{eq:Recurrence equation}
\ee
where $q_k$ are polynomials in 
$x$, and $m\in\N$.
If $q_0(n,k),q_m(n,k)\not\equiv 0$ then the number $m$ is called the order 
of the given recurrence equation.
\end{enumerate}
Some of those families, especially all ``classical'' families of orthogonal 
polynomials, have the following further property:
\begin{enumerate}
\item[4.]
{\bf (Rodrigues representation)}
\\
The functions $f_n$ have a representation of the Rodrigues type
\be
f_n(x)=\frac{1}{K_n\,g(x)}\dedn {x}{n} h_n(x)
\label{eq:Rodrigues}
\ee
for some functions $g$ depending on $x$, and $h_n$ depending on $n$ and $x$,
and a constant $K_n$ depending on $n$.
\end{enumerate}
From an algebraic point of view these properties read as follows:
Let $K[x]$ denote the field of rational functions over $K$ where $K$ is one
of $\Q$, $\R$, or $\C$. Then if the coefficients of the occurring
polynomials and rational functions are elements of $K$,
\begin{enumerate}
\item
the derivative rule states that $f_n'$ is an element of the linear space 
over $K[x]$ which is generated by $\{ f_n, f_{n-1}, \ldots, f_{n-(m-1)}\}$ or
$\{ f_n, f_{n+1}, \ldots, f_{n+{m-1}}\}$, respectively;
\item
the differential equation states that the $m+1$ functions $f_n^{(k)}\;
(k=0,\ldots,m)$ are linearly dependent over $K[x]$; moreover, by an induction
argument, any $m+1$ functions $f_n^{(k)}\;(k\in\N_0)$ are linearly dependent 
over $K[x]$;
\item
the recurrence equation states that the $m+1$ functions $f_{n-k}\;
(k=0,\ldots,m)$ are linearly dependent over $K[x]$; moreover, by an induction
argument, any $m+1$ functions $f_{n}\;(n\in D)$, are linearly dependent over 
$K[x]$.
\end{enumerate}
One important question when dealing with special functions is the
following: Which properties of those functions does one have to know 
to be able to establish various types of identities that those functions
satisfy? With respect to
the implementation of special functions in Computer Algebra systems this
question reads: Which properties should be implemented for those functions,
and in which form should this be done such that the user is enabled
to verify various types of identities, or at least to implement algorithms
for this purpose?

Nikiforov and Uvarov \cite{NU} 
gave a unified introduction to
special functions of mathematical physics based primarily on the Rodrigues
formula and the differential equation. They dealt, however, only with
second order differential equations, which makes their treatment quite
restricted, and moreover their development does not have algorithmic 
applications. 

Truesdell \cite{Truesdell} gave a unified approach to special functions
based entirely on a special form of the derivative rule.
His development has some algorithmic content, which, however, is difficult
or impossible to implement in Computer Algebra.
Truesdell's approach---although nice---has the further disadvantage that one 
can obtain only 
results of a very special form, see \cite{Koe94Truesdell}.

From the algorithmic point of view another approach is better:
We will base our treatment of special functions on the derivative
rule (\ref{eq:Derivative rule}) in combination with the recurrence equation
(\ref{eq:Recurrence equation}). We will show that an implementation
of special functions in Computer Algebra systems based on these two properties 
gives a simplification mechanism at hand which, in particular, enables
the user to verify many kinds of identities for those functions. Some of
these identities like differential equations, and power series
representations can even be dealt with algorithmically, i.\ e.\ they can
be computed by the Computer Algebra system.

Our treatment is connected with the holonomic system approach due to Zeilberger
\cite{Zei1}--\cite{Zei3} which is based on the valididy of partial 
differential equations, mixed recurrence equations, and
difference-differential equations. This connection will be made more
precise later.

The class of functions that can be treated this way contains the 
Airy functions $\airyai\:(x)$, $\airybi\:(x)$ (see e.\ g.\ \cite{AS}, \S~10.4), 
the Bessel functions $J_n(x), Y_n(x), I_n(x),$ and $K_n(x)$ 
(see e.\ g.\ \cite{AS}, Ch.~9--11), 
the Hankel functions $H_n^{(1)}(x)$ and $H_n^{(2)}(x)$
(see e.\ g.\ \cite{AS}, Ch.~9),
the Kummer functions $M(a,b,x)=\,_1 F_1\left.\lk\begin{array}{c}
\multicolumn{1}{c}{a}\\[-1mm] \multicolumn{1}{c}{b}
            \end{array}\rb x \rk$ and $U(a,b,x)$
(see e.\ g.\ \cite{AS}, Ch.~13),
the Whittaker functions $M_{n,m}(x)$ and $W_{n,m}(x)$
(see e.\ g.\ \cite{AS}, \S~13.4),
the associated Legendre functions $P_a^b(x)$ and $Q_a^b(x)$
(see e.\ g.\ \cite{AS}, \S~8),
all kinds of orthogonal polynomials:
the Jacobi polynomials $P_n^{(\al,\beta)}(x)$,
the Gegenbauer polynomials $C_n^{(\al)}(x)$,
the Chebyshev polynomials of the first kind $T_n(x)$ 
and of the second kind $U_n(x)$,
the Legendre polynomials $P_n(x)$,
the Laguerre polynomials $L_n^{(\al)}(x)$,
and the Hermite polynomials $H_n(x)$
(see \cite{Sze}, \cite{Tri}, and \cite{AS}, \S~22), 
many more special functions, and furthermore sums, products,
derivatives, antiderivatives,
and the composition with rational functions and rational powers
of those functions (see \cite{Sta}, \cite{Zei1}, \cite{SZ} and \cite{KS}).

In the case of the classical orthogonal polynomials the properties above can be 
made much more precise (see e.\ g.\ \cite{Tri}, Kapitel~IV).
Therefore let $f_n:[a,b]\pf\R\;(n\in\N_0)$ denote the family of orthogonal
polynomials 
\[
f_n(x)=k_n\,x^n+k'_n\,x^{n-1}+\ldots
\]
with respect to the weight function $w(x)\geq 0$, i.\ e.\ with the property
that
\[
\int\limits_a^b w(x)\,f_n(x)\,f_m(x)\,dx=0\quad\quad(n\neq m)
\]
and 
\[
\int\limits_a^b w(x)\,f_n^2(x)\,dx=h_n\neq 0
\;.
\]
Then we have the properties:
\begin{enumerate}
\item
{\bf (Derivative rule)}
\\
The functions $f_n$ satisfy a derivative rule of the form
\[
X\,f_n'=\beta_n\,f_{n-1}+\lk \frac{n}{2}X''x+\al_n\rk\,f_n
\]
(see e.\ g.\ \cite{Tri}, p.\ 135, formula (4.8))
where
\[
\al_n=n\,X'(0)-\ed 2\,X''\,\frac{k_n'}{k_n}
\;,\quad\quad
\beta_n=-\frac{h_n\,k_{n-1}}{h_{n-1}\,k_n}\lk
K_1\,k_1-\frac{2n-1}{2}\,X''\rk
\;,
\]
and
\be
X(x)=\funkdefff{(b-x)(x-a)}{a,b \;\mbox{are finite}}
{x-a}{b=\infty}{1}{-a,b=\infty}
\label{eq:X(x)}
\;.
\ee
Especially is the order of the derivative rule $2$.
\item
{\bf (Differential equation)}
\\
The functions $f_n$ satisfy the homogeneous
linear differential equation with polynomial coefficients
\[
X\,f_{n}''(x)+K_1\,f_1\,f_n'(x)+\la_n\,f_{n}(x)=0
\]
(see e.\ g.\ \cite{Tri}, p.\ 133, formula (4.1))
where
\[
\la_n=-n\lk K_1\,k_1-\frac{n-1}{2}\,X''\rk
\;,
\]
and $X(x)$ is given by (\ref{eq:X(x)}).
Especially is the order of the differential equation $2$.
\item
{\bf (Recurrence equation)}
\\
The functions $f_n$ satisfy the recurrence equation 
\be
f_{n+1}(x)=-C_n\,f_{n-1}(x)+(A_n\,x+B_n)\,f_n(x)
\label{eq:Recurrence equation:orth}
\ee
(see e.\ g.\ \cite{Tri}, p.\ 126, formula (2.1))
with 
\[
A_n=\frac{k_{n+1}}{k_n}\;,
\quad\quad
B_n=\frac{k_{n+1}}{k_n}\lk \frac{k'_{n+1}}{k_{n+1}}-\frac{k'_{n}}{k_{n}}\rk
\;,\quad\quad\mbox{and}\quad\quad
C_n=\frac{k_{n+1}\,k_{n-1}\,h_n}{k_n^2\,h_{n-1}}
\;.
\]
Especially is the order of the recurrence equation $2$.
\item
{\bf (Rodrigues representation)}
\\
The functions $f_n$ have a representation of the Rodrigues type
\be
f_n(x)=\frac{1}{K_n\,w(x)}\dedn {x}{n} \Big( w(x)\,X(x)^n \Big)
\label{eq:Rodriguestype}
\ee
(see e.\ g.\ \cite{Tri}, p.\ 129, formula (3.2)),
where $X(x)$ is given by (\ref{eq:X(x)}), i.\ e.\ 
(\ref{eq:Rodrigues}) is valid with $g(x)=w(x)$, and $h_n(x)=w(x)\,X(x)^n$. 
Especially: The order of the polynomial $X(x)$ is $\leq 2$.
\end{enumerate}
Further it turns out that in the case of classical orthogonal polynomials
all coefficient functions of $f_{n-k}$ are rational also
with respect to the variable $n$, a fact that depends, however,
on the special normalizations that are used in these cases.

We mention that no system of orthogonal polynomials besides the
classical ones satisfies a Rodrigues representation of type 
(\ref{eq:Rodriguestype}) with a polynomial $X$ (see e.\ g.\ \cite{Tri},
Kapitel IV, \S 3).

We note that using the recurrence equation (\ref{eq:Recurrence
equation:orth}), which is valid also for non-classical orthogonal 
polynomials, or any recurrence equation of type 
(\ref{eq:Recurrence equation}) of order two (also called three-term
recursion), recursively,
each (backward or forward) derivative rule (\ref{eq:Derivative rule}) is 
equivalent to a derivative rule
\be
f_n'(x)=k(n,x)\,f_n(x)+l(n,x)\,f_{n+1}(x)
\label{eq:Derivative rule2}
\ee
($k,l$ rational functions with respect to $x$)
of order two.
In general, the order of the derivative rule can always be assumed to
be less than or equal to the order of the recurrence equation.
In some nice work \cite{Truesdell} Truesdell presented a treatment of special
functions entirely based on the functional equation (\ref{eq:Derivative rule2}).
He showed that this difference-differential equation
is independent of the differential equation
(\ref{eq:Differential equation}) and the recurrence equation 
(\ref{eq:Recurrence equation}), i.\ e.\ it does not imply the existence
of one of these.

In contrast to this work, our main notion is the
\\[1mm]
{\bf Definition (Admissible family of special functions)}
We call a family $f_n 
$ of special functions {\sl admissible} if
the functions $f_n$ satisfy a recurrence equation of type
(\ref{eq:Recurrence equation}) and a derivative rule of type
(\ref{eq:Derivative rule}).
We call the order of the recurrence equation the {\sl order} of the admissible
family $f_n$.
\hfill$\Box$
\\[2mm]
Note that the recurrence equation (\ref{eq:Recurrence equation})
together with $m$ initial functions $f_{n_0},f_{n_0+1},\ldots,f_{n_0+m-1}$ 
determine the functions $f_n\;(n\in D)$ uniquely.

So an admissible family of special functions (with given initial
functions) is overdetermined by its two defining properties, i.\ e.\
the recurrence equation and the derivative rule must be compatible.
This fact, however, gives our notion a considerable strength:
\bT
\label{th:mdimensionsl}
{\rm
For any admissible family $f_n$ of order $m$ the linear space 
$V_{f_n}$ over $K[x]$
of functions generated by the set of shifted derivatives
$\{ f_{n\pm k}^{(j)}\;|\;j, k\in\N_0\}$ is at most $m$-dimensional.
On the other hand, if the family
$\{ f_{n\pm k}^{(j)}\;|\;j, k\in\N_0\}$ spans an $m$-dimensional linear
space, then $f_n$ forms an admissible family of order $m$.
}
\eT
\pr
By the recurrence equation and an induction argument it follows that 
the linear space $V$ spanned by $\{ f_{n\pm k}\;|\;k\in\N_0\}$ 
is at most $m$-dimensional. Using the derivative rule, 
by a further induction it follows that the derivative of any order
$f_n^{(k)}\;(k\in\N_0)$ is an element of $V$. Therefore $V_{f_n}=V$.

If on the other hand for a family $f_n$ the set of derivatives
$\{ f_{n\pm k}^{(j)}\;|\;j, k\in\N_0\}$ is $m$-dimensional, then
the existence of a recurrence equation and a derivative rule of
order $m$ are obvious.
\eop
From the algebraic point of view this is the main reason for the importance
of admissible families: Any $m+1$ distinguished elements of $V_{f_n}$ are
linearly dependent, i.\ e.\ any arbitrary element of $V_{f_n}$ can be
represented by a linear combination (with respect to $K[x]$) of any
$m$ of the others. This is the algebraic background for the fact that
so many identities between the members and their derivatives
of an admissible family exist.

In particular we have
\bC
\label{cor:AF->DE}
{\rm
Any admissible family $f_n$ of order $m$ satisfies a simple
differential equation of order $m$.
\hfill$\Box$
}
\eC
In \S~\ref{sec:Algorithmic generation of differential equations} we
give an algorithm which, in particular, generates this differential equation
of $f_n$.

With regard to Zeilberger's approach Corollary~\ref{cor:AF->DE} can be
interpreted as follows: Any admissible family $f_n(x)$ forms a holonomic
system with respect to the two variables $n$, and $x$, whose
defining recurrence equation, and the differential equation corresponding
to Corollary~\ref{cor:AF->DE} together with the initial conditions
\be
f_0^{(k)}(0)\;,\quad\quad\mbox{and}\quad\quad
f_k(0)\quad\quad(k=0,\ldots,m-1)
\label{eq:holonomicIV}
\ee
yield the canonical holonomic representation of $f_n(x)$ (see \cite{Zei1},
Lemma 4.1).

On the other hand, not all holonomic systems $f_n(x)$ form admissible 
families so that our notion is stronger: Let $f_n(x):=\airyai\:(x)$ for all
$n\in\Z$, then obviously $f_n(x)$ is the holonomic system generated by
the equations
\[
f_n''(x)=x\,f_n(x)\;,
\quad\quad
f_{n+1}(x)=f_n(x)\;,
\]
and some initial values,
that does {\sl not} form an admissible family as the derivative
$f_n'$ is linearly independent
of $\{f_n\;|\;n\in\Z\}$ over $K[x]$, see 
\S~\ref{sec:Embedding of one-variable functions into admissible families},
and thus no derivative rule of the form (\ref{eq:Derivative rule})
exists.

A further advantage of our approach is the separation of the variables,
i.\ e.\ the work with ordinary differential equations, and one-variable
recurrence equations rather than partial
differential equations, mixed recurrence equations, and
difference-differential equations.
So our approach---if applicable---seems to
be more natural.

To present an example of an admissible family
that cannot be found in mathematical dictionaries, we consider the functions
\[
k_n(x):=\frac{2}{\pi}\int\limits_0^{\pi/2}
\cos\:(x\,\tan\th-n\,\th)\,d\th
\;,
\]
that Bateman introduced in \cite{Bat}, see also \cite{KS1}.
He verified that (\cite{Bat}, formula (2.7))
\be
F_n(x):=(-1)^n\,k_{2n}(x)=(-1)^n\,e^{-x}\Big( L_n(2x)-L_{n-1}(2x)\Big)
\label{eq:k2m}
\;.
\ee
We call $F_n$ the family of Bateman functions which turns out
to be an admissible family of order two.

Bateman obtained the property (\cite{Bat}, formula (4.1))
\[
(n-1)\,\Big(F_n(x)-F_{n-1}(x)\Big)+(n+1)\,\Big(F_n(x)-F_{n+1}(x)\Big)=
2\,x\,F_n(x)
\]
leading to
\be
n\,F_{n}(x)-2\,(n-1-x)\,F_{n-1}(x)+(n-2)\,F_{n-2}(x)=0
\label{eq:differenceeq}
\ee
which is a recurrence equation of type (\ref{eq:Recurrence equation}) 
and order two that determines the Bateman functions uniquely using the
two initial functions
\[
F_0(x)=e^{-x}
\quad\quad\quad\mbox{and}\quad\quad\quad
F_1(x)=-2\,x\,e^{-x}
\]
which follow from (\ref{eq:k2m}). 

Bateman obtained further
a difference differential equation (\cite{Bat}, formula (4.2))
\be
(n+1)\,F_{n+1}(x)-(n-1)\,F_{n-1}(x)=2\,x\,F_n'(x)
\;,
\label{eq:difference differential equation}
\ee
which can be brought into the form
\be
F_n'(x)=\ed{x}\Big( (n-x)\,F_n(x)-(n-1)\,F_{n-1}(x)\Big)
\label{eq:Batemanderivativerule}
\ee
using (\ref{eq:differenceeq}). This is a derivative rule of the
form (\ref{eq:Derivative rule}) and order two. Therefore $F_n(x)$ form
an admissible family of order two.

We note that the functions $F_n$ satisfy the differential equation
\be
x\,F_n''(x)+(2n-x)\,F_n(x)=0
\;,
\label{eq:DE Bateman}
\ee
(see \cite{Bat}, formula (5.1)), and the Rodrigues type representation
\be
F_n(x)=\frac{x\,e^x}{n!}\frac{d^n}{dx^n}\lk e^{-2x}\,x^{n-1}\rk
\;,
\label{eq:Rodrigues Bateman}
\ee
(see \cite{Bat}, formula (31)).

\section{Properties of admissible families}
\label{sec:Properties of admissible families}
\bT
\label{th:Properties of admissible families}
{\rm
Let $f_n$ form an admissible family of order $m$. Then
\begin{enumerate}
\item[(a)]
{\bf (Shift)} 
$f_{n\pm k}\;(k\in\N)$ forms an admissible family of order $m$;
\item[(b)]
{\bf (Derivative)}
$f_n'$ forms an admissible family of order $\leq m$;
\item[(c)]
{\bf (Composition)}
$f_n\circ r$ forms an admissible family of order $\leq m$, 
if $r$ is a rational function, and of order $\leq m\,q$,
if $r(x)=x^{p/q}\;(p,q\in\N)$.
\end{enumerate}
If furthermore $g_n$ forms an admissible family of order $\leq l$, 
then moreover
\begin{enumerate}
\item[(d)]
{\bf (Sum)}
$f_n+g_n$ forms an admissible family of order $\leq m+l$;
\item[(e)]
{\bf (Product)}
$f_n\,g_n$ forms an admissible family of order $\leq m\,l$.
\end{enumerate}
}
\eT
\pr
(a): 
This is an obvious consequence of Theorem~\ref{th:mdimensionsl}.
\\(b):
Let $g_n:=f_n'$. We start with the recurrence equation for $f_n$
and take derivative to get
\be
\sum_{k=0}^m q_k'(n,x)\,f_{n-k}(x)+\sum_{k=0}^m q_k(n,x)\,f_{n-k}'(x)=0
\;.
\label{eq:intermediatederivative}
\ee
From Theorem~\ref{th:mdimensionsl}, we know that each of the functions
$f_{n-j}\;(j=0,\ldots,m)$ can be represented as a linear combination
of the functions $f_{n-k}'\;(k=0,\ldots,m-1)$ over $K[x]$, which
generates a recurrence equation for $g_n$. Similarly a derivative rule
for $g_n$ is obtained.
\\(c):
For the composition $h_n:=f_n\circ r$
with a rational function $r$, the recurrence equation
is obtained by substitution, and the derivative rule
is a result of the chain rule. If $r(x)=x^{1/q}$, then, by 
\cite{KS}, Lemma~1, the family $\{h_n^{(j)}\;|\;j\in\N_0\}$
is spanned by the $mq$ functions 
$x^{r/q} f_n^{(j)}(x^{1/q})\;(j=1,\ldots,m-1,$ $r=0,\ldots,q-1)$, and since
$\{f_{n\pm k}^{(j)}\;|\;j,k\in\N_0\}$ has dimension $m$, the linear
space spanned by
$\{h_{n\pm k}^{(j)}\;|\;j,k\in\N_0\}$ has dimension $\leq m\,q$,
implying the result. If finally $r(x)=x^{p/q}$, then a combination gives the 
result.
\\(d):
By a simple algebraic argument, we see that $f_{n-k}+g_{n-k}\:(k\in\Z)$ span
the linear space $V:=V_{f_n+g_n}=V_{f_n}+V_{g_n}$ of dimension $\leq m+l$
over $K[x]$.
Therefore $f_n+g_n$ satisfies a recurrence equation of order $\leq m+l$.
If we add the derivative rules for $f_n$ and $g_n$, we see that 
$f_n'+g_n'\in V$, and thus can be represented in the desired way.
\\(e):
By a similar algebraic argument (see e.\ g.\ \cite{Sta}, Theorem 2.3)
we see that $f_{n-k}\cdot g_{n-k}\:(k\in\Z)$ span a linear space $V$
of dimension $\leq m\,l$ over $K[x]$, hence $f_n\,g_n$ satisfies
a recurrence equation of order $\leq m\,l$. By the product rule,
and the derivative rules for $f_n$ and $g_n$
we see that the derivative of $f_n\,g_n$ is represented
by products of the form $f_{n-k}\,g_{n-j}\;(k,j\in\Z)$, and as 
those span the linear space $V$ (see e.\ g.\ \cite{KS}, Theorem 3 (d)), 
we are done.
\eop
As an application we again may state that the Bateman functions
form an admissible family: Using the theorem, this follows immediately from
representation (\ref{eq:k2m}).

Next we study algorithmic versions of the theorem. 
The following algorithm generates a representation of the members 
$f_{n\pm k}\;(k=0,\ldots,m-1)$ 
of an admissible family in terms of the derivatives 
$f_{n\pm j}'\;(j=0,\ldots,m-1)$. 
By Theorem~\ref{th:mdimensionsl} we know that such a representation exists.
Without loss of generality, we assume that the
admissible family is given by a backward derivative rule. In case of a
forward derivative rule, a similar algorithm is valid.

\begin{algorithm}
\label{algo:Integral rule}
{\rm
Let $f_n$ be an admissible family of order $m$, 
given by a backward derivative rule
\[
f_n'(x)=\sum_{k=0}^{m-1} r_k(n,x)\,f_{n-k}(x)
\;.
\]
Then the following algorithm generates a list of
backward rules $(k=0,\ldots,m-1)$
\be
f_{n-k}(x)=\sum_{j=0}^{m-1} R_j^k(n,x)\,f'_{n-j}(x)
\label{eq:Integral rule}
\ee
($R_j^k$ rational with respect to $x$) for $f_{n-k}\;(k=0,\ldots,m-1)$ 
in terms of the derivatives $f'_{n-j}\;(j=0,\ldots,m-1)$:
\begin{enumerate}
\item[(a)]
Shift the derivative rule $m-1$ times to obtain the set of $m$ equations
\[
f_{n-j}'(x)=\sum_{k=0}^{m-1} r_k(n-j,x)\,f_{n-j-k}(x)
\quad\quad(j=0,\ldots,m-1)
\;.
\]
\item[(b)]
Utilize the recurrence equation to express all expressions on the right hand
sides of these equations in terms of $f_{n-k}\;(k=0,\ldots,m-1)$ leading to
\[
f_{n-j}'(x)=\sum_{k=0}^{m-1} r_k^{j}(n,x)\,f_{n-k}(x)
\quad\quad(j=0,\ldots,m-1\;,\;r_k^{j}\;\mbox{rational 
with respect to}\;x)
\;.
\]
\item[(c)]
Solve this linear equations system for the variables $f_{n-k}\;(k=0,\ldots,m-1)$
to obtain the representations (\ref{eq:Integral rule}) searched for.
\hfill$\Box$
\end{enumerate}
}
\end{algorithm}
The proof of the algorithm is obvious. It is also clear how the method
can be adapted to obtain forward rules in terms of the derivatives.
As an example, the algorithm generates the following
representations for the Bateman functions
\[
F_n(x)=
  {{1 - n + x}\over {2\,n -1- x}}\,F_n'(x) +
   {{n - 1}\over {2\,n -1- x}}\,F_{n-1}'(x)
\;,
\]
and
\[
F_n(x)=
  {{1 + n - x}\over {1 + 2\,n - x}}\,F_n'(x) -
   {{1 + n}\over {1 + 2\,n - x}}\,F_{n+1}'(x)
\]
in terms of their derivatives.

We note that by means of Algorithm~\ref{algo:Integral rule}
and the results of \cite{KS} (see also \cite{Zei1}, p.\ 342, and \cite{SZ}),
we are able to state algorithmic versions of the statements of
Theorem~\ref{th:Properties of admissible families}.

\begin{algorithm}
\label{algo:Properties of admissible families}
{\rm
The following algorithms lead to the derivative rules
and recurrence equations of the admissible families presented in
Theorem~\ref{th:Properties of admissible families}:
\begin{enumerate}
\item[(a)]
{\bf (Shift)}
Direct use of derivative rule and recurrence equation lead to the 
derivative rule and the recurrence equation for $f_{n\pm 1}$; a recursive
application gives the results for $f_{n\pm k}\;(k\in\N)$.
\item[(b)]
{\bf (Derivative)}
By Algorithm~\ref{algo:Integral rule} we may  replace all occurrences
of $f_{n-k}\;(k=0,\ldots,m)$ in (\ref{eq:intermediatederivative}),
resulting in the recurrence equation for $f_n'$; similarly the derivative
rule is obtained.
\item[(c)]
{\bf (Composition)}
If $r$ is a rational function, then an application of the chain rule leads
to the derivative rule and the recurrence equation of $f_n\circ r$;
an approach similar to the algorithmic version of Theorem 2 in \cite{KS}
yields the derivative rule and the recurrence equation of $f_n\circ x^{1/q}$
by an elimination of the expressions $x^{r/q}\,f_n^{(j)}(x^{1/q})\;
(r=1,\ldots,q-1,\;j=1,\ldots,m-1)$.
\item[(d)]
{\bf (Sum)}
Applying a discrete version of Theorem 3 (c) in \cite{KS}
to $f_n+g_n$ (see also \cite{Zei1}, p.\ 342, and
\cite{SZ}, {\sc Maple} function {\tt rec+rec})
results in the recurrence equation, and a similar approach gives 
the derivative rule.
\item[(e)]
{\bf (Product)}
Applying a discrete version of Theorem 3 (d) in \cite{KS}
to $f_n\,g_n$ (see also \cite{Zei1}, p.\ 342, and 
\cite{SZ}, {\sc Maple} function {\tt rec*rec})
yields the recurrence equation, and a similar approach gives
the derivative rule.
\hfill$\Box$
\end{enumerate}
}
\end{algorithm}
A {\sc Mathematica} implementation of the given algorithms generate
e.\ g.\ for the derivative $F_n'(x)$ of the Bateman function $F_n(x)$ 
the derivative rule
\[
F_n''(x)= \frac{2\,n - x}{x - 2\,n\,x + x^2}\Big(
\lk n-1  \rk F_{n-1}'(x)+ \lk 1-n+x  \rk F_{n}'(x)
\Big)
\;,
\]
and the recurrence equation
\[
F_{n+1}'(x)=
\ed{(1\! + \!n) (1\! -\! 2 n\! + \!x)}
\lk
(n\!-\!1)(x\!-\!2n\!-\!1) F_{n-1}'(x)+ 
2\,(1\!-\!2n^2\!+\! 3nx\!-\!x^2) F_{n}'(x)
\rk
\;,
\]
and for the product $A_n(x):=F_n^2(x)$ the derivative rule
\beao
A_n'(x)&=&
{{\left( 1 - n \right) \,{{\left( n-2 \right) }^2}}\over 
{2\,n\,x\,\left( 1 - n + x \right) }}\,A_{n-2}(x) 
\\&&+\; 
{{2\,\left( n-1 \right) \,\left( 1 - n + x \right)}\over {n\,x}} 
\,A_{n-1}(x)
\\&&+ \;
{{\left( 3\,n - 3\,{n^2} - 4\,x + 8\,n\,x - 4\,{x^2} \right)}\over 
{2\,x\,\left( 1 - n + x \right) }}\,A_n(x)
\;,
\eeao
and the recurrence equation
\beao
A_{n+1}(x) &= &
\ed{(1 + n)^2}\,\lk
\frac{(n-2)^2\,(n-1)\,(x-n)}{n\,(1 - n + x)}\,A_{n-2} 
\right. 
\\&&+\;
\frac{(n-1)\,(3\,n - 3\,n^2 - 4\,x + 8\,n\,x - 4\,x^2)}{n}\,A_{n-1} 
\\&&+\;
\left.
\frac{(x-n)\,(-3\,n + 3\,n^2 + 4\,x - 8\,n\,x + 4\,x^2)}{1 - n + x}\,A_n
\rk
\eeao
are derived.

\section{Derivative rules of special functions}
\label{sec:Derivative rules of special functions}

Many Computer Algebra systems like {\sc Axiom} \cite{Axi},
{\sc Macsyma} \cite{Mac}, {\sc Maple} \cite{Map}, {\sc Mathematica}
\cite{Wol}, or {\sc Reduce} \cite{Red} support the work with
special functions. On the other hand, there are so many identities for
special functions that it is a nontrivial task to decide which
properties should be used by the system (and in which way)
for the work with those functions. 

Since all Computer Algebra systems support derivatives, as a first question
it is natural to 
ask how the current implementations of Computer Algebra systems handle
the derivatives of special functions. Here are some examples:
{\sc Mathematica} (Version 2.2) gives

{\small
\begin{verbatim}
In[1]:= D[BesselI[n,x],x]

        BesselI[-1 + n, x] + BesselI[1 + n, x]
Out[1]= --------------------------------------
                          2
In[2]:= D[LaguerreL[n,a,x],x]

Out[2]= -LaguerreL[-1 + n, 1 + a, x]
\end{verbatim}
}\noindent
We note that in {\sc Mathematica} the derivatives of all special
functions symbolically are implemented. 
On the other hand, we notice that, given the function $I_n\:(x)$,
{\sc Mathematica}'s derivative introduces two new functions: $I_{n-1}\:(x)$, 
and $I_{n+1}\:(x)$. Given the Laguerre polynomial $L_n^{(\al)}(x)$, the
derivative produced introduces a new function where both $n$, and $\al$
are altered. The representation used is optimal for numerical purposes, 
but is not a representation according to our classification.

With {\sc Maple} (Version V.2) we get

{\small
\begin{verbatim}
> diff(BesselI(n,x),x);
                                          n BesselI(n, x)
                      BesselI(n + 1, x) + ---------------
                                                 x

> diff(L(n,a,x),x);
                                  d
                                ---- L(n, a, x)
                                 dx
\end{verbatim}
}\noindent
Thus {\sc Maple}'s derivative for the Bessel function $I_n\:(x)$ 
introduces only one new function $I_{n+1}(x)$, and is
of type (\ref{eq:Derivative rule}), whereas 
(even if {\tt orthopoly} is loaded) no symbolic derivative of 
the Laguerre polynomial $L_n^{(\al)}(x)$ is implemented.

Obviously there is no unique way to declare the derivative of a special 
function. However, we note that if we declare the derivative of a
special function by a derivative rule of type
(\ref{eq:Derivative rule}) of order $m$ then we can be
sure that the derivative
of the special function $f_n(x)$ introduces at most $m$ new functions,
namely $f_{n-k}(x)\;(k=1,\ldots,m)$. Moreover, if the family of
special functions depends on several parameters, then the given representation 
of the derivative does not use any functions with other parameters changed.


Here we give a list of the backward derivative rules of the form 
(\ref{eq:Derivative rule}) for the families of special functions that
we introduced in \S~\ref{sec:Introduction} which all turn out to be
of order two (see 
e.\ g.\ \cite{AS}, 
(9.1.27) (Bessel and Hankel functions), 
(9.2.26) (Bessel functions), 
(13.4.11), (13.4.26) (Kummer functions), 
(13.4.29)--(13.4.33) (Whittaker functions),
(8.5.4) (associated Legendre functions), 
and \S~22.8 (orthogonal polynomials)):
%
\bea
J_n'\:(x)
&=&
J_{n-1}\:(x)-\frac{n}{x}\,J_n\:(x)
\;,
\label{eq:Jnstrich}
\nonumber
\\[3mm]
Y_n'\:(x)
&=&
Y_{n-1}\:(x)-\frac{n}{x}\,Y_n\:(x)
\;,
\label{eq:Ynstrich}
\nonumber
\\[3mm]
I_n'\:(x)
&=&
I_{n-1}\:(x)-\frac{n}{x}\,I_n\:(x)
\;,
\label{eq:Instrich}
\nonumber
\\[3mm]
K_n'\:(x)
&=&
-K_{n-1}\:(x)-\frac{n}{x}\,K_n\:(x)
\;,
\label{eq:Knstrich}
\nonumber
\\[3mm]
\ded x H_n^{(1)}(x)
&=&
H_{n-1}^{(1)}\:(x)-\frac{n}{x}\,H_n^{(1)}\:(x)
\label{eq:Hn1strich}
\;,
\nonumber
\\[3mm]
\ded x H_n^{(2)}(x)
&=&
H_{n-1}^{(2)}\:(x)-\frac{n}{x}\,H_n^{(2)}\:(x)
\label{eq:Hn2strich}
\;,
\nonumber
\\[3mm]
\ded x M(a,b,x)
\!\!&=&
\ed x\Big( (b-a)\,M(a-1,b,x)-(b-a-x)\,M(a,b,x) \Big)
\label{eq:Mstrich}
\;,
\nonumber
\\[3mm]
\ded x U(a,b,x)
\!\!&=&
\ed x\Big( -U(a-1,b,x)+(a-b+x)\,U(a,b,x)\Big)
\label{eq:KummerUstrich}
\;,
\nonumber
\\[3mm]
%
%
%
%
%
%
%
%
%
%
%
%
%
M_{n,m}'\:(x)
&=&
\frac{1}{2x}\Big( (1+2m-2n)\,M_{n-1,m}\:(x)+(2n-x)\,M_{n,m}\:(x)\Big)
\;,
\label{eq:Mnmstrich}
\nonumber
\\[3mm]
%
%
%
%
%
%
%
%
%
%
%
%
%
%
W_{n,m}'\:(x)
&=&
\ed{4x}\lk (1-4m^2-4n+4n^2)\,W_{n-1,m}\:(x)+(4n-2x)\,W_{n,m}\:(x)\rk
\;,
\label{eq:Wnmstrich}
\nonumber
\\[3mm]
\ded x P_a^b(x)
&=&
\ed{1-x^2}\lk (a+b)\,P_{a-1}^b(x)-a\,x\,P_a^b(x)\rk
\;,
\label{eq:Pabstrich}
\nonumber
\\[3mm]
\ded x Q_a^b(x)
&=&
\ed{1-x^2}\lk (a+b)\,Q_{a-1}^b(x)-a\,x\,Q_a^b(x)\rk
\;,
\label{eq:Qabstrich}
\nonumber
\\[3mm]
\ded x P_n^{(\al,\beta)}(x)
\!\!
&=&
\!\!
\ed{(2n\!+\!\al\!+\!\beta)(1\!-\!x^2)}
\lk
2(n\!+\!\al)(n\!+\!\beta)P_{n-1}^{(\al,\beta)}(x)+
n(\al\!-\!\beta\!-\!(2n\!+\!\al\!+\!\beta)x)P_{n}^{(\al,\beta)}(x)
\rk
\!,
\label{eq:Jacobistrich}
\nonumber
\\[3mm]
\ded x C_n^{(\al)}\:(x)
&=&
\ed{1-x^2}
\lk
(n+2\al-1)\,C_{n-1}^{(\al)}\:(x)-n\,x\,C_n^{(\al)}\:(x)
\rk
\;,
\label{eq:Gegenbauerstrich}
\nonumber
\\[3mm]
T_n'\:(x)
&=&
\ed{1-x^2}
\Big(
n\,T_{n-1}\:(x)-n\,x\,T_n\:(x)
\Big)
\;,
\label{eq:ChebyshevTstrich}
\nonumber
\\[3mm]
U_n'\:(x)
&=&
\ed{1-x^2}
\Big(
(n+1)\,U_{n-1}\:(x)-n\,x\,U_n\:(x)
\Big)
\;,
\label{eq:ChebyshevUstrich}
\nonumber
\\[3mm]
P_n'\:(x)
&=&
\ed{1-x^2}
\Big(
n\,P_{n-1}\:(x)-n\,x\,P_n\:(x)
\Big)
\;,
\label{eq:Legendrestrich}
\nonumber
\\[3mm]
\ded x L_{n}^{(\al)}(x)
&=&
\ed{x}
\lk
-(n+\al)\,L_{n-1}^{(\al)}(x)+n\,L_{n}^{(\al)}(x)
\rk
\;,
\label{eq:Laguerrestrich}
\\[3mm]
H_{n}'(x)
&=&
2n\,H_{n-1}(x)
\;.
\label{eq:Hermitestrich}
\nonumber
\eea

\section{Recurrence equations of special functions}
\label{sec:Recurrence equations of special functions}

Whenever in any expression 
subexpressions of the form $r_k\,f_{n-k}\;(r_k\;\mbox{rational}, k\in\Z)$ 
occur, in an admissible family of order $m$ 
with the recursive use of the recurrence equation we may replace
so many occurrences of those expressions $r_k\,f_{n-k}$
that finally only $m$ successive terms of the same type remain.

This allows for example to eliminate the number of occurrences
in any linear combination (over $K[x]$)
of derivatives of $f_n$ to $m$, a fact with which we will deal in more detail
in \S~\ref{sec:Algorithmic generation of differential equations}.

We show how {\sc Mathematica} and {\sc Maple} work with regard to this
question. Whereas {\sc Mathematica} does not have any built-in
capabilities to simplify the following linear combinations of Bessel and
Laguerre functions,
{\small
\begin{verbatim}
In[3]:= BesselI[n+1,x]+2*n/x*BesselI[n,x]-BesselI[n-1,x]

                              2 n BesselI[n, x]
Out[3]= -BesselI[-1 + n, x] + ----------------- + BesselI[1 + n, x]
                                      x

In[4]:= Simplify[%]
                              2 n BesselI[n, x]
Out[4]= -BesselI[-1 + n, x] + ----------------- + BesselI[1 + n, x]
                                      x

In[5]:= LaguerreL[n+1,a,x]-(2*n+a+1-x)*LaguerreL[n,a,x]+(n+a)*LaguerreL[n-1,a,x]

Out[5]= (a + n) LaguerreL[-1 + n, a, x] - 
 
>    (1 + a + 2 n - x) LaguerreL[n, a, x] + LaguerreL[1 + n, a, x]

In[6]:= Simplify[%]

Out[6]= (a + n) LaguerreL[-1 + n, a, x] - 
 
>    (1 + a + 2 n - x) LaguerreL[n, a, x] + LaguerreL[1 + n, a, x]
\end{verbatim}
}\noindent
with {\sc Maple} we get

{\small
\begin{verbatim}
> BesselI(n+1,x)+2*n/x*BesselI(n,x)-BesselI(n-1,x);

                                n BesselI(n, x)
          BesselI(n + 1, x) + 2 --------------- - BesselI(- 1 + n, x)
                                       x

> simplify(");
                                       0

> L(n+1,a,x)-(2*n+a+1-x)*L(n,a,x)+(n+a)*L(n-1,a,x);

     L(n + 1, a, x) - (2 n + a + 1 - x) L(n, a, x) + (n + a) L(n - 1, a, x)

> simplify(");

   L(n + 1, a, x) - 2 L(n, a, x) n - L(n, a, x) a - L(n, a, x) + L(n, a, x) x

        + L(n - 1, a, x) n + L(n - 1, a, x) a
\end{verbatim}
}\noindent
%
i.\ e.\ {\sc Maple}'s {\tt simplify} command supports simplification with the 
aid of the recurrence equations for the Bessel functions.
On the other hand, for the orthogonal polynomials (even if {\tt orthopoly} 
is loaded) no simplifications occur.

In the rest of this section we give a list of the recurrence equations of 
the given type for the families of special functions that we consider 
which all turn out to be of order two (see e.\ g.\ \cite{AS},
(9.1.27), 
(9.2.26), 
(13.4.1), (13.4.15), 
(13.4.29), (13.4.31), 
(8.5.3), 
and \S~22.7). 
We list them in the form explicitly solved for $F_{n+1}$ as this is the usual 
form found in mathematical dictionaries.
\bea
J_{n+1}\:(x)
&=&
-J_{n-1}\:(x)+\frac{2n}{x}\,J_n\:(x)
\;,
\label{eq:Jn+1}
\nonumber
\\[3mm]
Y_{n+1}\:(x)
&=&
-Y_{n-1}\:(x)+\frac{2n}{x}\,Y_n\:(x)
\;,
\label{eq:Yn+1}
\nonumber
\\[3mm]
I_{n+1}\:(x)
&=&
I_{n-1}\:(x)-\frac{2n}{x}\,I_n\:(x)
\;,
\label{eq:In+1}
\nonumber
\\[3mm]
K_{n+1}\:(x)
&=&
K_{n-1}\:(x)+\frac{2n}{x}\,K_n\:(x)
\;,
\label{eq:Kn+1}
\nonumber
\\[3mm]
H_{n+1}^{(1)}\:(x)
&=&
-H_{n-1}^{(1)}\:(x)+\frac{2n}{x}\,H_n^{(1)}\:(x)
\;,
\label{eq:Hn+11}
\nonumber
\\[3mm]
H_{n+1}^{(2)}\:(x)
&=&
-H_{n-1}^{(2)}\:(x)+\frac{2n}{x}\,H_n^{(2)}\:(x)
\;,
\label{eq:Hn+12}
\nonumber
\\[3mm]
M(a+1,b,x)
&=&
\ed a\Big( (b-a)\,M(a-1,b,x)+(2a-b+x)\,M(a,b,x)
\Big)
\;,
\label{eq:Mn+1}
\nonumber
\\[3mm]
U(a+1,b,x)
&=&
-\frac{1}{a\,(1+a-b)}\Big(
U(a-1,b,x)+(b-2a-x)\,U(a,b,x)
\Big)
\;,
\label{eq:Un+1}
\nonumber
\\[3mm]
%
%
%
%
%
%
%
M_{n+1,m}\:(x)
&=&
\frac{1}{1 + 2 m + 2 n}
\Big(
(1+2m-2n)\,M_{n-1,m}\:(x)+(4n-2x)\,M_{n,m}\:(x)
\Big)
\;,
\label{eq:Mn+1m}
\nonumber
\\[3mm]
%
%
%
%
%
%
%
W_{n+1,m}\:(x)
&=&
\ed{4}
\lk 
(-1+4m^2+4n-4n^2)\,W_{n-1,m}\:(x)-(8n-4x)\,W_{n,m}\:(x)
\rk
\;,
\label{eq:Wn+1m}
\nonumber
\\[3mm]
P_{a+1}^b\:(x)
&=&
\ed{a-b+1}
\Big(
-(a+b)\,P_{a-1}^b\:(x)+(2a+1)\,x\,\,P_{a-1}^b\:(x)
\Big)
\;,
\label{eq:Pa+1b}
\nonumber
\\[3mm]
Q_{a+1}^b\:(x)
&=&
\ed{a-b+1}
\Big(
-(a+b)\,Q_{a-1}^b\:(x)+(2a+1)\,x\,\,Q_{a-1}^b\:(x)
\Big)
\;,
\label{eq:Qn+1b}
\nonumber
\\[3mm]
P_{n+1}^{(\al,\beta)}(x)
&=&
\ed{2\,(n\!+\!1)\,(n\!+\!\al\!+\!\beta\!+\!1)\,(2n\!+\!\al\!+\!\beta)}
\lk
-2(n\!+\!\al)(n\!+\!\beta)(2n\!+\!\al\!+\!\beta\!+2)\,P_{n-1}^{(\al,\beta)}(x)
\right.
\nonumber
\\[3mm]
&&
+
\left.
\lk (2n\!+\!\al\!+\!\beta\!+1)(\al^2-\beta^2)+
(2n\!+\!\al\!+\!\beta)_3\,x\rk P_{n}^{(\al,\beta)}(x)
\rk
\;,
\label{eq:Jacobin+1}
\nonumber
\\[3mm]
C_{n+1}^{(\al)}\:(x)
&=&
\ed{n+1}
\lk
-(n+2\al-1)\,C_{n-1}^{(\al)}\:(x)+2(n+\al)\,x\,C_n^{(\al)}\:(x)
\rk
\;,
\label{eq:Gegenbauern+1}
\nonumber
\\[3mm]
T_{n+1}\:(x)
&=&
-T_{n-1}\:(x)+2\,x\,T_n\:(x)
\;,
\label{eq:ChebyshevTn+1}
\nonumber
\\[3mm]
U_{n+1}\:(x)
&=&
-U_{n-1}\:(x)+2\,x\,U_n\:(x)
\;,
\label{eq:ChebyshevUn+1}
\nonumber
\\[3mm]
P_{n+1}\:(x)
&=&
\ed{n+1}
\Big(
-n\,P_{n-1}\:(x)+(2n+1)\,x\,P_n\:(x)
\Big)
\;,
\label{eq:Legendren+1}
\nonumber
\\[3mm]
L_{n+1}^{(\al)}(x)
&=&
\ed{n+1}
\lk
-(n+\al)\,L_{n-1}^{(\al)}(x)+(2n+\al+1-x)\,L_{n}^{(\al)}(x)
\rk
\;,
\label{eq:Laguerren+1}
\nonumber
\\[3mm]
H_{n+1}(x)
&=&
-2n\,H_{n-1}(x)+2x\,H_{n}(x)
\;.
\label{eq:Hermiten+1}
\nonumber
\eea
Note that $(a)_k$ (which is used in the recurrence equation for the
Jacobi polynomials $P_{n}^{(\al,\beta)}(x)$) denotes the 
{\sl Pochhammer symbol\/} (or {\sl shifted factorial\/}) defined by
$(a)_{k}:=\prod\limits_{j=1}^k (a\!+\!j\!-\!1)$.

We note further that for functions with several ``discrete'' variables 
it may happen that for each of them there exists a recurrence equation.
As an example we consider the Laguerre polynomials for which we have
(\cite{AS} (22.7.29), in combination with (22.7.30))
\be
L_{n}^{(\al+1)}(x)
=
\ed{x}
\lk
-(n+\al)\,L_{n}^{(\al-1)}(x)+(\al+x)\,L_{n}^{(\al)}(x)
\rk
\;.
\label{eq:alternateRELaguerre}
\ee
In \S~\ref{sec:Functions of the hypergeometric type as admissible families}
we will demonstrate that generalized hypergeometric functions
satisfy recurrence equations with respect to all their parameters.

To be safely enabled that the algorithms of 
\S~\ref{sec:Algorithmic verification of identities}--%
\S~\ref{sec:Algorithmic verification of formulas involving symbolic sums}
apply, all of those 
recurrence equations should be implemented and applied recursively
for simplification purposes.

\section{Embedding of one-variable functions into admissible families}
\label{sec:Embedding of one-variable functions into admissible families}

In this section we consider first, how the elementary transcendental
functions are covered by the given approach.

Consider the exponential function $f(x)=e^x$. This function can be embedded
into the admissible family $f_n$, defined by the properties
\[
f_n'(x)=f_n(x)\;,
\quad\quad
f_{n+1}(x)=f_n(x)
\quad\quad\mbox{and}\quad\quad
f_0(x)=e^x
\;,
\]
i.\ e.\ the family of iterated derivatives of $e^x$.

Obviously this is a representation of an admissible family of order one.

Moreover in the given case it turns out that $f_n(x)=e^x=f_0(x)$ for all
$n\in\Z$, so there is no actual need to give the functions numbers, 
and therefore we (obviously) keep the usual notation.

Similarly the functions $\sin x$ and $\cos x$ are embedded into the admissible 
family $f_n$ of order two given by the properties
\[
f_n'(x)=f_{n-1}(x)\;,
\quad\quad
f_{n+1}(x)=-f_{n-1}(x)
\;,
\quad\quad\mbox{and}\quad\quad
f_0(x)=\cos x\;,
\quad\quad
f_1(x)=\sin x
\;.
\]
Again, the family of functions $f_n$ is finite, and our numbering is
unnecessary:
\[
f_n(x)=\funkdeffff
{\cos x}{n=4m\;(m\in\Z)}
{\sin x}{n=4m+1\;(m\in\Z)}
{-\cos x}{n=4m+2\;(m\in\Z)}
{-\sin x}{n=4m+3\;(m\in\Z)}
\;.
\]
Essentially there are only the two functions $\cos x$, and $\sin x$
involved. Note, however, that both functions are needed as
no simple first order differential equation for $\sin x$ or $\cos x$
exists.

Other nontrivial
examples of essentially finite admissible families of special functions
are formed by the Airy functions. Let $\airyai_n\:(x)=\airyai^{(n)}\:(x)$,
i.\ e.\
\[
\airyai_n'\:(x)=\airyai_{n+1}\:(x)
\;.
\]
By the differential equation for
the Airy functions (see e.\ g.\ \cite{AS}, (10.4)) we have
$\airyai''\:(x)-x\,\airyai\:(x)=0$, so that from Leibniz's rule
it follows that
\beao
\airyai_{n+1}\:(x)
&=&
\airyai^{(n+1)}\:(x)=
\Big( \airyai''\:(x)\Big)^{(n-1)}
\\&=&
\Big( x\,\airyai\:(x)\Big)^{(n-1)}=
\sum_{k=0}^{n-1} \ueber{n-1}{k}\,x^{(k)}\,\Big(\airyai\:(x)\Big)^{(n-1-k)}
\\&=&
x\,\airyai^{(n-1)}\:(x)+(n-1)\,\airyai^{(n-2)}\:(x)=
x\,\airyai_{n-1}\:(x)+(n-1)\,\airyai_{n-2}\:(x)
\;,
\eeao
and therefore $\airyai\:(x)$ 
is embedded into the admissible family $\airyai_n$ of order three given by
\be
\airyai_n'\:(x)=\airyai_{n+1}\:(x)\;,
\quad\quad
\airyai_{n+1}\:(x)=x\,\airyai_{n-1}\:(x)+(n-1)\,\airyai_{n-2}\:(x)
\;,
\label{eq:DR,RE AiryAi}
\ee
and we have the initial functions
\[
\airyai_0\:(x)=\airyai\:(x)\;,
\quad\quad
\airyai_1\:(x)=\airyai'\:(x)
\quad\quad\mbox{and}\quad\quad
\airyai_2\:(x)=x\,\airyai\:(x)
\;.
\]
Similarly $\airybi\:(x)$ is embedded into the admissible family 
of order three given by
\be
\airybi_n'\:(x)=\airybi_{n+1}\:(x)\;,
\quad\quad
\airybi_{n+1}\:(x)=x\,\airybi_{n-1}\:(x)+(n-1)\,\airybi_{n-2}\:(x)
\;,
\label{eq:DR,RE AiryBi}
\ee
and the initial functions
\[
\airybi_0\:(x)=\airybi\:(x)\;,
\quad\quad
\airybi_1\:(x)=\airybi'\:(x)
\quad\quad\mbox{and}\quad\quad
\airybi_2\:(x)=x\,\airybi\:(x)
\;.
\]
Our indexed families turn out to be representable by
\[
\airyai_n\:(x)=p_n(x)\,\airyai\:(x)+q_n(x)\,\airyai'\:(x)
\quad\quad\mbox{and}\quad\quad
\airybi_n\:(x)=p_n(x)\,\airybi\:(x)+q_n(x)\,\airybi'\:(x)
\;,
\]
with polynomials $p_n$ and $q_n$ in $x$.
This shows, however, that to deal with the Airy functions
algorithmically as is suggested in this paper, besides the functions
$\airyai\:(x)$ and $\airybi\:(x)$ the two {\sl independent} functions
$\airyai'\:(x)$ and $\airybi'\:(x)$ are needed, but none else. 
Let's look, how Computer Algebra systems work with the Airy functions.

{\sc Maple} handles them as follows:

{\small
\begin{verbatim}
> Ai(x);
                                     Ai(x)
> diff(Ai(x),x);

       1/2                   3/2
      2    BesselK(1/3, 2/3 x   )
  1/4 ---------------------------
                 1/4
                x    Pi

                       /                                                 3/2 \
              1/2  5/4 |                    3/2        BesselK(1/3, 2/3 x   )|
             2    x    |- BesselK(4/3, 2/3 x   ) + 1/2 ----------------------|
                       |                                         3/2         |
                       \                                        x            /
       + 1/3 -----------------------------------------------------------------
                                             Pi
\end{verbatim}
}\noindent

{\small
\begin{verbatim}
> simplify(diff(Ai(x),x$2)-x*Ai(x));

             1/2                   3/2       1/2                    3/2   3/2
  1/48 (- 3 2    BesselK(1/3, 2/3 x   ) - 8 2    BesselK(-2/3, 2/3 x   ) x

             1/2  3                   3/2        9/4             /   5/4
       + 16 2    x  BesselK(1/3, 2/3 x   ) - 48 x    Ai(x) Pi)  /  (x    Pi)
                                                               /
> diff(Bi(x),x);
                                     d
                                   ---- Bi(x)
                                    dx
> diff(Bi(x),x$2);
                                     2
                                    d
                                  ----- Bi(x)
                                     2
                                   dx
\end{verbatim}
}\noindent
So the derivative of $\airyai\:(x)$ is represented by Bessel functions,
whereas the function $\airyai\:(x)$ itself is not, and therefore the
expression {\tt diff(Ai(x),x\verb+$+2)-x*Ai(x)} is not simplified. 
On the other hand the derivative of $\airybi\:(x)$ is {\sl not} a 
valid {\sc Maple} function. With {\sc Mathematica} we get
 
{\small
\begin{verbatim}
In[7]:= D[AiryAi[x],x]

Out[7]= AiryAiPrime[x]

In[8]:= D[AiryAiPrime[x],x]

Out[8]= x AiryAi[x]

In[9]:= D[AiryAi[x],{x,2}]-x*AiryAi[x]

Out[9]= 0

In[10]:= D[AiryBi[x],x]

Out[10]= AiryBiPrime[x]

In[11]:= D[AiryBiPrime[x],x]

Out[11]= x AiryBi[x]

In[12]:= D[AiryBi[x],{x,2}]-x*AiryBi[x]

Out[12]= 0
\end{verbatim}
}\noindent
Thus we see that in this situation
{\sc Mathematica} does exactly what we suggest: It works with the independent
functions $\airyai\:(x)$, $\airyai'\:(x)$, $\airybi\:(x)$, $\airybi'\:(x)$,
and the derivative rules (\ref{eq:DR,RE AiryAi}) and (\ref{eq:DR,RE AiryBi}).

As a further example of an admissible family we consider the iterated
integrals
\[
\erfc_n\:(x)=\int\limits_x^\infty \erfc_{n-1}\:(t)\,dt
\]
of the (complementary) error function $\erfc\:(x)=1-\erf\:(x)=\erfc_0\:(x)$
(see e.\ g.\ \cite{AS}, (7.2)) that form the admissible family with
\[
\erfc_n'\:(x)=-\erfc_{n-1}\:(x)\;,
\quad\quad
\erfc_{n+1}\:(x)=\frac{1}{2(n+1)}\erfc_{n-1}\:(x)-\frac{x}{n+1}\,\erfc_{n}\:(x)
\;,
\]
and the initial functions
\[
\erfc_0\:(x)=\erfc\:(x)\;,
\quad\quad
\erfc_1\:(x)=-\ed{\sqrt\pi} \left(\sqrt\pi\,x\,\erfc\:(x)-e^{-x^2}\right)
\]
(one may also use the initial value function
$\erfc_{-1}\;(x)=\frac{2}{\sqrt\pi}e^{-x^2}$).
In particular, $\erfc x$ is embedded into an admissible family.

{\sc Maple} deals with these functions as suggested:

{\small
\begin{verbatim}
> diff(erfc(n,x),x);
                                - erfc(n - 1, x)

> simplify(diff(erfc(n,x),x$2)+2*x*diff(erfc(n,x),x)-2*n*erfc(n,x));

                                       0
\end{verbatim}
}\noindent
As a final example, we mention another family of iterated integrals,
the {\sl Abramowitz functions}
\[
A_n(x):=\int\limits_0^\infty t^n\,e^{-t^2-x/t}\,dt
\]
(see \cite{Abr}, and \cite{AS}, (27.5)) 
which form an admissible family with derivative rule
\[
A_n'(x)=
\ded x \lk \int\limits_0^\infty t^n\,e^{-t^2-x/t}\,dt\rk
=
\int\limits_0^\infty \ded x \lk t^n\,e^{-t^2-x/t}\rk dt
=
-\int\limits_0^\infty t^{n-1}\,e^{-t^2-x/t}\,dt=
-A_{n-1}(x)
\]
of order one (see \cite{AS}, (27.5.2)), and recurrence formula
\[
A_{n+1}(x)=\frac{n}{2}\,A_{n-1}(x)+\frac{x}{2}\,A_{n-2}(x)
\]
of order three (\cite{AS}, (27.5.3)).

Again, embedded into an admissible family, especially the function
$A_0(x)=\int\limits_0^\infty e^{-t^2-x/t}\,dt$ is covered by our approach.

%

\section{Embedding the inhomogeneous case}
\label{sec:Embedding the inhomogeneous case}

Some families of functions are characterized by inhomogeneous
differential rules and recurrence equations. 
Examples for this situation are the exponential integrals
given by
\[
E_n\:(x)=\int\limits_1^\infty \frac{e^{xt}}{t^n}\,dt
\]
(see e.\ g.\ \cite{AS}, (5.1)), and the Struve functions
${\bf H}_n(x)$ and ${\bf L}_n(x)$ (see e.\ g.\ \cite{AS}, Chapter~5),
for which we have the inhomogeneous properties
\[
E_n'\:(x)=-E_{n-1}\:(x)\;,
\quad\quad
E_{n+1}\:(x)=\frac{e^{-x}}{n}-\frac{x}{n}\,E_n\:(x)
\;,
\]
(\cite{AS}, (5.1.14) and (5.1.26)),
\be
{\bf H}_{n-1}(x)-{\bf H}_{n+1}(x)=2\,{\bf H}_n'(x)-
\frac{x^n}{2^n\,\sqrt\pi\:\Gamma(n+3/2)}
\;,
\label{eq:Hnstrichorig}
\ee
\[
{\bf H}_{n-1}(x)+{\bf H}_{n+1}(x)=\frac{2n}{x}\,{\bf H}_n(x)+
\frac{x^n}{2^n\,\sqrt\pi\:\Gamma(n+3/2)}
\]
(\cite{AS}, (12.1.9)--(12.1.10)),
and
\be
{\bf L}_{n-1}(x)+{\bf L}_{n+1}(x)=2\,{\bf L}_n'(x)-
\frac{x^n}{2^n\,\sqrt\pi\:\Gamma(n+3/2)}
\;,
\label{eq:Lnstrichorig}
\ee
\[
{\bf L}_{n-1}(x)-{\bf L}_{n+1}(x)=\frac{2n}{x}\,{\bf L}_n(x)+
\frac{x^n}{2^n\,\sqrt\pi\:\Gamma(n+3/2)}
\]
(\cite{AS}, (12.2.4)--(12.2.5)), respectively.
Eliminating the inhomogeneous parts (using 
$\Gamma(3/2+n)=(1/2+n)\,\Gamma(1/2+n)$), 
these examples are made into admissible
families with the derivative rules
\bea
E_n'\:(x)
&=&
-E_{n-1}\:(x)
\;,
\label{eq:Einstrich}
\nonumber
\\[3mm]
{\bf H}_{n}'(x)
&=&
{\bf H}_{n-1}(x)-\frac{n}{x}\,{\bf H}_{n}(x)
\;,
\label{eq:StruveHstrich}
\\[3mm]
{\bf L}_{n}'(x)
&=&
{\bf L}_{n-1}(x)-\frac{n}{x}\,{\bf L}_{n}(x)
\;,
\label{eq:StruveLstrich}
\eea
and the recurrence equations
\bea
E_{n+1}\:(x)
&=&
\ed{n}\Big( x\,E_{n-1}(x) +(n-1-x)\,E_n(x) \Big)
\;,
\label{eq:Einrecurrence}
\nonumber
\\[3mm]
{\bf H}_{n+1}(x)
&=&
\ed{2n+1}\Big( 
x\,{\bf H}_{n-2}(x)+(1-4n)\,{\bf H}_{n-1}(x)+
\frac{x^2+2 n + 4 n^2}{x}\,{\bf H}_{n}(x)
\Big)
\;,
\label{eq:StruveHrecurrence}
\nonumber
\\[3mm]
{\bf L}_{n+1}(x)
&=&
\ed{2n+1}\Big(
-x\,{\bf L}_{n-2}(x)-(1-4n)\,{\bf L}_{n-1}(x)+
\frac{x^2-2 n-4 n^2}{x}\,{\bf L}_{n}(x)
\Big)
\;,
\label{eq:StruveLrecurrence}
\nonumber
\eea
so that the exponential integrals form an admissible family of order
two, and the Struve functions ${\bf H}_{n}(x)$ and ${\bf L}_{n}(x)$
form admissible families of order three. Note that the above derivative
rules (\ref{eq:StruveHstrich})--(\ref{eq:StruveLstrich})
are not listed in \cite{AS} although they are much simpler than
the inhomogeneous relations 
(\ref{eq:Hnstrichorig})--(\ref{eq:Lnstrichorig}).

After bringing the inhomogeneous rules into the desired form, those families
are recognized as admissible families, and our method can be applied.

\section{Functions of the hypergeometric type as admissible families}
\label{sec:Functions of the hypergeometric type as admissible families}

All functions introduced in this paper are special cases of functions of the
hypergeometric type (see \cite{Koe92}). In this section we will show that the
generalized hypergeometric function $_{p}F_{q}$ defined by
\be
_{p}F_{q}\left.\lk\begin{array}{cccc}
a_{1}&a_{2}&\cdots&a_{p}\\
b_{1}&b_{2}&\cdots&b_{q}\\
            \end{array}\rb x\rk
:=
\sumi A_k\,x^{k}=
\sumi \frac
{(a_{1})_{k}\cdot(a_{2})_{k}\cdots(a_{p})_{k}}
{(b_{1})_{k}\cdot(b_{2})_{k}\cdots(b_{q})_{k}\,k!}x^{k}
\label{eq:coefficientformula}
\;,
\ee
and thus by Theorem~\ref{th:Properties of admissible families}~(c)
all functions of the hypergeometric type, form admissible 
families. Therefore we first deduce a derivative rule of order two
for $_{p}F_{q}$.

Let us choose any of the 
numerator parameters $n:=a_k\;(k=1,\ldots,p)$ of $_{p}F_{q}$ as parameter $n$.
Further we use the abbreviations
\[
F_n(x)=\;_{p}F_{q}\left.\lk\begin{array}{cccc}
n &a_{2}&\cdots&a_{p}\\
b_{1}&b_{2}&\cdots&b_{q}\\
            \end{array}\rb x\rk
=
\sumi A_k(n)\,x^{k}
\;.
\]
From the relation
\[
\frac{(n+1)_k}{(n)_k}=\frac{n+k}{n}
\]
it follows that
\[
n\,A_{k}(n+1)=(n+k)\,A_{k}(n)
\;.
\]
Using the differential operator $\theta f(x)=x\,f'(x)$,
we get by summation
\beao
n\,F_{n+1}(x)&=&n\sumi A_{k}(n+1)\,x^{k}=
(n+k)\sumi A_{k}(n)\,x^{k}
\\&=&
n\,F_{n}(x)+\sumi k\,A_k(n)\,x^{k}=
n\,F_{n}(x)+\theta F_n(x)
\;,
\eeao
and therefore we are led to the derivative rule
\be
\theta F_n(x)=n\,\Big(F_{n+1}(x)-F_n(x)\Big)\;,
\quad\quad\mbox{or}\quad\quad
F_n'(x)=\frac{n}{x}\,\Big(F_{n+1}(x)-F_n(x)\Big)\;.
\label{eq:hypergeoDR}
\ee
Hence we have established that for any of the
numerator parameters $n:=a_k\;(k=1,\ldots,p)$ of $_{p}F_{q}$ such a
simple (forward) derivative rule is valid.

We note that by similar means for each of the denominator
parameters $n:=b_k\;(k=1,\ldots,q)$ of $_{p}F_{q}$ the simple (backward)
derivative rule 
\be
\theta F_n(x)=(n-1)\,\Big(F_{n-1}(x)-F_n(x)\Big)\;,
\quad\quad\mbox{or}\quad\quad
F_n'(x)=\frac{n-1}{x}\,\Big(F_{n-1}(x)-F_n(x)\Big)
\label{eq:hypergeoDRbackward}
\ee
is derived.

Next, we note that $F_n$
satisfies the well-known hypergeometric differential equation
\be
\th (\th+b_1-1)\cdots (\th+b_q-1)F_{n}(x)
=x(\th+a_1)(\th+a_2)\cdots(\th+a_p)F_{n}(x)
\;.
\label{eq:hypergeoDE}
\ee
Replacing all occurrences
of $\theta$ in (\ref{eq:hypergeoDE}) recursively by the derivative rule
(\ref{eq:hypergeoDR}) or (\ref{eq:hypergeoDRbackward}), 
a recurrence equation for $F_n$ is obtained
that turns out to have the same order as the differential
equation (\ref{eq:hypergeoDE}), i.\ e.\ $\max\{p,q+1\}$. 

We summarize the above results in the following
\bT
\label{th:generalized hypergeometric function}
{\rm
The generalized hypergeometric function 
$\;_{p}F_{q}\left.\lk\begin{array}{cccc}
a_{1}&a_{2}&\cdots&a_{p}\\
b_{1}&b_{2}&\cdots&b_{q}\\
            \end{array}\rb x\rk$
satisfies the derivative rules
\[
\theta F_n(x)=n\,\Big(F_{n+1}(x)-F_n(x)\Big)
\]
for any of its
numerator parameters $n:=a_k\;(k=1,\ldots,p)$, and
\[
\theta F_n(x)=(n-1)\,\Big(F_{n-1}(x)-F_n(x)\Big)
\]
for any of its denominator parameters $n:=b_k\;(k=1,\ldots,q)$,
and recursive substitution of all occurrences of $\theta$ in the hypergeometric
differential equation
\[
\th (\th+b_1-1)\cdots (\th+b_q-1)F_{n}(x)
=x(\th+a_1)(\th+a_2)\cdots(\th+a_p)F_{n}(x)
\]
generates a recurrence equation of the type (\ref{eq:Recurrence equation})
of order $\max\{p,q+1\}$ with respect to the parameter chosen.
This recurrence equation has coefficients that are rational
with respect to $x$, and $n$.
In particular, $_{p}F_{q}$ forms an admissible family 
of order $\max\{p,q+1\}$ with respect to all of its parameters
$a_k,b_k$.
\hfill$\Box$
}
\eT
We note that
if some of the parameters of $_{p}F_{q}$ are specified, there may exist
a lower order differential equation, and thus the order of the admissible
family may be lower than the theorem states. We note further 
that this theorem is the main reason for the fact that so many special
functions form admissible families: Most of them can be represented in
terms of generalized hypergeometric functions. 

\section{Algorithmic generation of differential equations}
\label{sec:Algorithmic generation of differential equations}

In this section we show that the algorithm to generate the uniquely
determined differential
equation of type (\ref{eq:Differential equation}) of 
lowest order valid for $f$ which was developed in \cite{Koe92} (see also
\cite{KS}), does apply if $f$ is constructed from functions that are
embedded into admissible families.

\begin{algorithm}[Find a simple differential equation]
\label{algo:Find a simple DE}%
{\rm
Let $f$ be a function given by an expression
that is built from the functions $\exp x$, $\ln x$, 
$\sin x$, $\cos x$, $\arcsin x$, $\arctan x$, and any other functions that
are embedded into admissible families, with the aid of the following 
procedures: differentiation, antidifferentiation,
addition, multiplication, and the composition with
rational functions and rational powers. 

Then the following procedure generates 
a simple differential equation valid for $f$:
\begin{enumerate}
\item[{\rm (a)}]
Find out whether there exists a simple differential equation
for $f$ of order $N:=1$.
Therefore differentiate $f$, and solve the linear equation
\[
f'(x)+A_{0}f(x)=0
\]
for $A_{0}$; i.\ e.\ set $A_{0}:=-\frac{f'(x)}{f(x)}$.
Is $A_{0}$ rational in $x$, then you are done after
multiplication with its denominator.
\item[{\rm (b)}]
Increase the order $N$ of the differential equation searched for by one.
Expand the expression
\[
f^{(N)}(x)+A_{N-1}f^{(N-1)}(x)+\cdots+A_{0}f(x)
\;,
\]
apply the recurrence formulas of any admissible family $F_{n}$
of order $m$ involved recursively
to minimize the occurrences of $F_{n-k}$ to at most $m$ successive
$k$-values, and check, if the remaining summands contain exactly $N$
rationally independent expressions considering
the numbers $A_{0}, A_{1},\ldots, A_{N-1}$ as constants.
Just in that case there exists a solution as follows: Sort with respect to
the rationally independent terms and create a system of linear equations
by setting their coefficients to zero. Solve this system for the numbers
$A_0, A_1,\ldots, A_{N-1}$. Those are rational functions in $x$, and
if there is a solution, this solution is unique. 
After multiplication by the common denominator of
$A_0, A_1,\ldots, A_{N-1}$ you get the differential equation searched for.
Finally cancel common factors of the polynomial coefficients.
\item[{\rm (c)}]
If part (b) was not successful, repeat step (b).
\end{enumerate}
}
\end{algorithm}
\pr
Theorem~3 of \cite{KS} (compare \cite{Sta}) shows that for $f$ a differential
equation of type (\ref{eq:Differential equation}) exists. 
We assume that differentiation is done by recursive descent through the
expression tree, and an application of the chain, product and quotient 
rules on the corresponding subexpressions. It is clear that the algorithm 
works for members of admissible families, compare Theorem~\ref{th:mdimensionsl}
and Corollary~\ref{cor:AF->DE}. Similarly the algorithm obviously
works for derivatives and antiderivatives of admissible families.
Further it is easily seen that the derivatives of sums, products,
and the composition with rational functions and rational powers form
either sums, or sums of products
all of which by a recursive use of the recurrence equations involved
are represented by sums of fixed lengths,
compare Theorem~\ref{th:Properties of admissible families}.
Thus after a finite number of steps, part (b) of the algorithm will succeed
(sharp a priory bounds for the resulting orders are given in \cite{KS}).
\eop
We note that from the implementational point of view the crucial step
of the algorithm is the decision of the rational independency in part (b).
If this decision can be handled properly, then the proof given in 
\cite{Koe92} shows that 
the algorithm generates the simple differential equation of lowest order 
valid for $f$. 

In our implementations, for testing whether some terms are rationally 
dependent, we divide each one by any other and test whether the quotient 
is a rational function in $x$ 
or not. This is an easy and fast approach which never leads
to wrong results, but may miss a simpler solution, which in practice, 
rarely happens.

Typically this happens, however, for orthogonal polynomials with prescribed 
$n$, for which a first order differential equation exists.
In this case, the recurrence equation hides these
rational dependencies, and in some sense (s.\ \cite{GK}, \S~7) 
here it is even advantageous that the rational dependency is not realized.

Another example where our implementations yield a differential equation
which is not of lowest order is given by
{\small
\begin{verbatim}
In[13]:= SimpleDE[Sin[2 x]-2 Sin[x] Cos[x],x]

Out[13]= 4 F[x] + F''[x] == 0
\end{verbatim}
}\noindent
This happens because the functions $\sin\:(2x)$ and $2\,\sin x\cos x$
algebraically cannot be verified to be rationally dependent even though they
are identical.

We note that, for elementary functions, we could use the Risch normalization
procedure~\cite{Risch} to generate the rationally independent
terms, but this does not work for special functions.

Further we note that in case of expressions of high complexity, the
use of \cite{KS}, Algorithm 2, typically is faster.
This algorithm, however, in general leads to a
differential equation of higher order than 
Algorithm~\ref{algo:Find a simple DE}.

As a first application of Algorithm~\ref{algo:Find a simple DE}
we consider the Airy functions $\airyai_n$, again, for
which the {\sc Mathematica} implementation of our algorithm yields

{\small
\begin{verbatim}
In[14]:= SimpleDE[AiryAi[n,x],x]

                                    (3)
Out[14]= (-1 - n) F[x] - x F'[x] + F   [x] == 0
\end{verbatim}
}\noindent
i.\ e.\ the differential equation
\be
\airyai_{n}'''\:(x)-x\,\airyai_{n}'\:(x)-(n+1)\,\airyai_{n}\:(x)=0
\;.
\label{eq:DE AiryAi}
\ee
Similarly, we get for the square of the Airy function

{\small
\begin{verbatim}
In[15]:= SimpleDE[AiryAi[x]^2,x]

                                (3)
Out[15]= -2 F[x] - 4 x F'[x] + F   [x] == 0
\end{verbatim}
}\noindent
The next calculation confirms the differential equation for the Bateman
functions $F_n$ (\ref{eq:DE Bateman})

{\small
\begin{verbatim}
In[16]:= SimpleDE[Bateman[n,x],x]

Out[16]= (2 n - x) F[x] + x F''[x] == 0
\end{verbatim}
}\noindent
Other examples are given with the aid of the iterated
integrals of the complementary error function, and the Abramowitz functions:

{\small
\begin{verbatim}
In[17]:= SimpleDE[Erfc[n,x],x]

Out[17]= -2 n F[x] + 2 x F'[x] + F''[x] == 0
\end{verbatim}
}\noindent
(see \cite{AS} (7.2.2))
and

{\small
\begin{verbatim}
In[18]:= SimpleDE[Exp[a x]*Erfc[n,x],x]

           2
Out[18]= (a  - 2 n - 2 a x) F[x] + (-2 a + 2 x) F'[x] + F''[x] == 0

In[19]:= SimpleDE[Exp[a x^2]*Erfc[n,x],x]

                            2      2  2
Out[19]= (-2 a - 2 n - 4 a x  + 4 a  x ) F[x] + (2 x - 4 a x) F'[x] +

>     F''[x] == 0

In[20]:= SimpleDE[Abramowitz[n,x],x]

                                      (3)
Out[20]= 2 F[x] + (1 - n) F''[x] + x F   [x] == 0
\end{verbatim}
}\noindent
(see \cite{AS} (26.2.41)).

We note that the algorithm obviously works for antiderivatives.
An example of that type is Dawson's integral (see e.\ g.\ \cite{AS}
(7.1.17)) for which we get the differential equation

{\small
\begin{verbatim}
In[21]:= SimpleDE[E^(-x^2)*Integrate[E^(t^2),{t,0,x}],x]

Out[21]= 2 F[x] + 2 x F'[x] + F''[x] == 0
\end{verbatim}
}\noindent
For the Struve functions, our algorithm generates the differential equations
%
%
%
\[
(n^2+n^3+x^2-nx^2)\,{\bf H}_n(x)+x\,(x^2-n-n^2)\,{\bf H}_n'(x)+
(2-n)\,x^2\,{\bf H}_n''(x)+x^3\,{\bf H}_n'''(x)
=0
\;,
\]
%
and 
\[
(n^2+n^3-x^2+nx^2)\,{\bf L}_n(x)-x\,(x^2+n+n^2)\,{\bf L}_n'(x)+
(2-n)\,x^2\,{\bf L}_n''(x)+x^3\,{\bf L}_n'''(x)
=0
\;,
\]
that are the homogeneous counterparts of the differential equation
(12.1.1) in \cite{AS}.

Finally we give examples involving hypergeometric functions:

{\small
\begin{verbatim}
In[22]:= SimpleDE[Hypergeometric2F1[a,b,c,x],x]

Out[22]= a b F[x] + (-c + x + a x + b x) F'[x] + (-1 + x) x F''[x] == 0

In[23]:= SimpleDE[Hypergeometric2F1[a,b,a+b+1/2,x]^2,x]

Out[23]= 8 a b (a + b) F[x] + 2 (-a - 2 a  - b - 4 a b - 2 b  + x + 3 a x +

            2                          2
>        2 a  x + 3 b x + 8 a b x + 2 b  x) F'[x] +

>     3 x (-1 - 2 a - 2 b + 2 x + 2 a x + 2 b x) F''[x] +

                  2  (3)
>     2 (-1 + x) x  F   [x] == 0
\end{verbatim}
}\noindent
Here the last function considered
$
\lk
_{2}F_{1}\left.\lk\begin{array}{c}
\multicolumn{1}{c}{\begin{array}{cc}a&b\end{array}}\\
\multicolumn{1}{c}{a\!+\!b\!+\!1/2}\\
            \end{array}\rb x\rk
\rk^2
$
is the left hand side of Clausen's formula (\ref{eq:Clausen's formula})
that we will consider again in
\S~\ref{sec:Algorithmic verification of identities}.

Now we investigate the case that a derivative rule and a differential equation
are given, and show that these two imply the existence of a recurrence equation:
\begin{algorithm}
{\rm
If a family $f_n$ is given by a derivative rule of type
(\ref{eq:Derivative rule}) and a differential equation of type 
(\ref{eq:Differential equation}), then it forms an admissible family
for which a recurrence equation can be found algorithmically.
}
\end{algorithm}
\pr
We present an algorithm which generates a recurrence equation for $f_n$:
Iterative differentiation of the derivative rule (\ref{eq:Derivative rule})
with the explicit use of (\ref{eq:Derivative rule}) at each step yields
\[
f_n^{(j)}(x)=\sum_{k=0}^M r_k^{j}(n,x)\,f_{n-k}(x)
\]
with rational functions $r_k^{j}$.
The substitution of these derivative representations in the differential 
equation gives the recurrence equation searched for.
\eop
As an example we consider the Airy functions $\airyai_n$, again, for
which we have the derivative rule (\ref{eq:DR,RE AiryAi})
\[
\airyai_n'\:(x)=\airyai_{n+1}\:(x)
\]
and the differential equation (\ref{eq:DE AiryAi})
\[
\airyai_{n}'''\:(x)-x\,\airyai_{n}'\:(x)-(n+1)\,\airyai_{n}\:(x)=0
\;.
\]
Differentiating the derivative rule successively and substituting the 
resulting expressions into the differential equation immediately yields the
recurrence equation (\ref{eq:DR,RE AiryAi}), again.

If this family, however, is given by the backward derivative rule 
(compare (\ref{eq:DR,RE AiryAi}))
\[
\airyai_n'\:(x)=x\,\airyai_{n-1}\:(x)+(n-1)\,\airyai_{n-2}\:(x)\;,
\]
then differentiation yields
\beao
\airyai_n''\:(x)&=&\airyai_{n-1}\:(x)+x\,\airyai_{n-1}'\:(x)+
(n-1)\,\airyai_{n-2}'\:(x)
\\&=&
\airyai_{n-1}\:(x)\!+\!x\Big( x\airyai_{n-2}\:(x)\!
+\!(n\!-\!2)\airyai_{n-3}\:(x)
\Big)
\!+\!(n\!-\!1)\Big( x\airyai_{n-3}\:(x)\!+\!(n\!-\!3)\airyai_{n-4}\:(x)\Big)
\\&=&
\airyai_{n-1}\:(x)+
x^2\,\airyai_{n-2}\:(x)+
(2n - 3 )\,x\,\airyai_{n-3}\:(x)+
(n^2-4 n+3)\,\airyai_{n-4}\:(x)
\;.
\eeao
After a similar procedure we get
\beao
\airyai_n'''\:(x)&=&
3x\,\airyai_{n-2}\:(x)+
(x^2+3n-5)\,\airyai_{n-3}\:(x)+
(3n-6)\,x^2\,\airyai_{n-4}\:(x)
\\&&+
(3n^2-15n+15)\,x\,\airyai_{n-5}\:(x)+
(n^3-9n^2+23n-15)\,\airyai_{n-6}\:(x)
\;,
\eeao
and the substitution into the differentiation equation gives finally
\beao
&&(n-1)\,\airyai_{n}\:(x)
-x^2\,\airyai_{n-1}\:(x)+
(4-n)\,x\,\airyai_{n-2}\:(x)+
(3n-5+x^3)\,\airyai_{n-3}\:(x)
\\&&+
(3n\!-\!6)\,x^2\,\airyai_{n-4}\:(x)+
(3n^2\!-\!15n\!+\!15)\,x\,\airyai_{n-5}\:(x)+
(n^3\!-\!9n^2\!+\!23n\!-\!15)\,\airyai_{n-6}\:(x)
=0
\;,
\eeao
a recurrence equation of order 6 rather than the minimal order three.
This shows, that, in general, the order of the resulting recurrence equation
is not best possible.

Algebraically spoken, our result tells that if $\{ f_n^{(j)}\;|\;j\in\N_0\}$
has finite dimension, and if $f_n'$ is an element of the linear space $V$
spanned by a finite number of the functions $\{f_{n\pm k}\}$,
then the space generated by all of $\{f_{n\pm k}\}$ is
of finite dimension, too. In contrast to Theorem~\ref{th:mdimensionsl}, however,
the dimension of this space generally may be higher than the dimension of $V$.
This shows the advantage of the use of admissible families.

As a further result of this section
we note that using our general procedure developed in \cite{Koe92} we have 
\begin{algorithm}[Find a Laurent-Puiseux representation]
\label{algo:Find a Laurent-Puiseux representation}%
{\rm
Let $f$ be a function that is built from the functions $\exp x$, $\ln x$,
$\sin x$, $\cos x$, $\arcsin x$, $\arctan x$, and any other functions that
are embedded into admissible families, with the aid of the following
procedures: differentiation, anti\-differentiation, addition,
multiplication, and the composition with
rational functions and rational powers.

If furthermore $f$ turns out to be of rational, exp-like, or
hypergeometric type (see \cite{Koe92}), then a closed form 
Laurent-Puiseux representation $f(x)=\sum\limits_{k=k_0}^\infty a_k\,x^{k/n}$
can be obtained algorithmically.
\hfill$\Box$
}
\end{algorithm}
We remark that there is a decision procedure due to Petkovsek \cite{P}
to decide the hypergeometric type from the recurrence equation obtained.

With Algorithm~\ref{algo:Find a Laurent-Puiseux representation},
it is possible to reproduce most of the results of the
extensive bibliography on series \cite{Han}, and to generate others.
As an example we present the power series representation of the square
of the Airy function:

{\small
\begin{verbatim}
In[24]:= PowerSeries[AiryAi[x]^2,x]

                1 k   k  1 + 3 k
               (-)  27  x        (2 k)!
                9
Out[24]= Sum[-(------------------------), {k, 0, Infinity}] +
               Sqrt[3] Pi k! (1 + 3 k)!

           k  3 k            1
         12  x    Pochhammer[-, k]
                             6
>    Sum[-------------------------, {k, 0, Infinity}] +
             1/3              2 2
          3 3    (3 k)! Gamma[-]
                              3

            1/3   k          2 + 3 k            5
         2 3    12  (1 + k) x        Pochhammer[-, k]
                                                6
>    Sum[--------------------------------------------, {k, 0, Infinity}]
                                      1 2
                     (3 + 3 k)! Gamma[-]
                                      3
\end{verbatim}
}\noindent
Note that, moreover, this technique generates hypergeometric 
representations, whenever such representations exist. The above example,
e.\ g., is recognized as the hypergeometric representation
\beao
\airyai\:(x)^2
&=&
\ed{3^{4/3}\,\Gamma\:(2/3)^2}\;
_{1}F_{2}\left.\lk\begin{array}{c}
\multicolumn{1}{c}{ 1/6 }\\[1mm]
\multicolumn{1}{c}{\begin{array}{cc}1/3 & 2/3\end{array}}\\
            \end{array}\rb \frac{4}{9}\,x^3 \rk
\\&&
-\frac{x}{\sqrt 3\,\pi}\;
_{1}F_{2}\left.\lk\begin{array}{c}
\multicolumn{1}{c}{1/2}\\[1mm]
\multicolumn{1}{c}{\begin{array}{cc}2/3 & 4/3 \end{array}}\\
            \end{array}\rb \frac{4}{9}\,x^3 \rk
+
\frac{x^2}{3^{2/3}\,\Gamma\:(1/3)^2}\;
_{1}F_{2}\left.\lk\begin{array}{c}
\multicolumn{1}{c}{ 5/6 }\\[1mm]
\multicolumn{1}{c}{\begin{array}{cc}4/3 & 5/3\end{array}}\\
            \end{array}\rb \frac{4}{9}\,x^3 \rk
\;.
\eeao
As soon as a hypergeometric representation is obtained, by
Theorem~\ref{th:generalized hypergeometric function} derivatives rules
and recurrence equations with respect to all parameters involved may be
obtained. As an example, we consider the Laguerre polynomials: The power series
representation for the Laguerre polynomial $L_n^{(\al)}(x)$ 
that our algorithm generates corresponds to the hypergeometric 
representation
\[
L_n^{(\al)}(x)=\ueber{n+\al}{n}\;
_1 F_1\left.\lk\begin{array}{c}
\multicolumn{1}{c}{-n}\\[1mm]
\multicolumn{1}{c}{\al+1}\\
            \end{array}\rb x \rk
\]
from which by an application of 
Theorem~\ref{th:generalized hypergeometric function}
we obtain the derivative rule
\beao
\ded x L_n^{(\al)}(x)
&=&
\ueber{n+\al}{n}\,\frac{-n}{x}\lk
_1 F_1\left.\lk\begin{array}{c}
\multicolumn{1}{c}{-n+1}\\[1mm]
\multicolumn{1}{c}{\al+1}\\
            \end{array}\rb x \rk
-\;
_1 F_1\left.\lk\begin{array}{c}
\multicolumn{1}{c}{-n}\\[1mm]
\multicolumn{1}{c}{\al+1}\\
            \end{array}\rb x \rk
\rk
\\&=&
\frac{-(n+\al)}{x}\,\ueber{n-1+\al}{n-1}\;
_1 F_1\left.\lk\begin{array}{c}
\multicolumn{1}{c}{-(n-1)}\\[1mm]
\multicolumn{1}{c}{\al+1}\\
            \end{array}\rb x \rk
+\frac{n}{x}\,\ueber{n+\al}{n}\;
_1 F_1\left.\lk\begin{array}{c}
\multicolumn{1}{c}{-n}\\[1mm]
\multicolumn{1}{c}{\al+1}\\
            \end{array}\rb x \rk
\\&=&
\ed{x}\lk -(n+\al)\,L_{n-1}^{(\al)}(x)+n\,L_n^{(\al)}(x)\rk
\;,
\eeao
i.\ e.\ (\ref{eq:Laguerrestrich}), again, but we are also led to
the derivative rule with respect $\al$:
\beao
\ded x L_n^{(\al)}(x)
&=&
\ueber{n+\al}{n}\,\frac{\al}{x}\lk
_1 F_1\left.\lk\begin{array}{c}
\multicolumn{1}{c}{-n}\\[1mm]
\multicolumn{1}{c}{\al}\\
            \end{array}\rb x \rk
-\;
_1 F_1\left.\lk\begin{array}{c}
\multicolumn{1}{c}{-n}\\[1mm]
\multicolumn{1}{c}{\al+1}\\
            \end{array}\rb x \rk
\rk
\\&=&
\frac{\al}{x}\,\frac{n+\al}{\al}\,\ueber{n+\al-1}{n}\,
_1 F_1\left.\lk\begin{array}{c}
\multicolumn{1}{c}{-n}\\[1mm]
\multicolumn{1}{c}{\al}\\
            \end{array}\rb x \rk
-
\frac{\al}{x}\,\ueber{n+\al}{n}\;
_1 F_1\left.\lk\begin{array}{c}
\multicolumn{1}{c}{-n}\\[1mm]
\multicolumn{1}{c}{\al+1}\\
            \end{array}\rb x \rk
\\&=&
\ed x\lk (n+\al)\,L_n^{(\al-1)}(x)-\al\,L_n^{(\al)}(x)\rk
\;.
\eeao
A further application of Theorem~\ref{th:generalized hypergeometric function}
yields the recurrence equation 
\[
F_{\al+1} = 
\frac{1+\al}{(1 + \alpha + n)\,x}\Big(
-\alpha\,F_{\alpha-1}+(\alpha + x)\,F_\alpha
\Big)
\]
for
$F_\al:=\,_1 F_1\left.\lk\begin{array}{c}
\multicolumn{1}{c}{-n}\\[-1mm] \multicolumn{1}{c}{\al+1}
            \end{array}\rb x \rk$
with respect to $\alpha$, and the use of the algorithm for the product
(\cite{KS}, Theorem 3 (d), \cite{Zei1}, p.\ 342, 
and \cite{SZ}, {\sc Maple} function {\tt rec*rec}),
applied to $L_n^{(\al)}=\ueber{n+\al}{n}\cdot F_\al$ generates 
(\ref{eq:alternateRELaguerre}), again.

\section{Algorithmic verification of identities}
\label{sec:Algorithmic verification of identities}

On the lines of \cite{Zei1} we can now present
an implementable algorithm to verify identities between expressions
using the results of the last section.

\begin{algorithm}
\label{algo:Verification of identities}
{\bf (Verification of identities)}
{\rm
Assume two functions $f_n(x)$ and $g_n(x)$ are given, to which 
Algorithm~\ref{algo:Find a simple DE} applies. Then the following
procedure verifies whether $f_n$ and $g_n$ are identical:
\begin{enumerate}
\item[(a)]
{\tt de1:=SimpleDE(f,x)}: \\
Determine the simple differential equation {\tt de1} corresponding to $f_n$.
\item[(b)]
{\tt de2:=SimpleDE(g,x)}: \\
Determine the simple differential equation {\tt de2} corresponding to $g_n$.
\item[(c)]
{\bf (Different differential equation implies different function)}
If {\tt de1} and {\tt de2} have the same order, then
\begin{itemize}
\item[-]
if they do not coincide besides common factors, i.\ e.\ have rational ratio,
then $f_n$ and $g_n$ do not coincide; return this, and quit.
\item[-]
Otherwise $f_n$ and $g_n$ satisfy the same differential equation {\tt de1} 
of order $l$, say, and it remains to check $l$ initial values. 
Continue with (e).
\end{itemize}
\item[(d)]
Let the orders of {\tt de1} and {\tt de2}, i.\ e.\
\[
\sum_{j=0}^l p_j\,f_n^{(j)}=0
\quad\quad\mbox{and}\quad\quad
\sum_{k=0}^{m} q_k\, g_n^{(k)}=0
\]
($p_j\;(j=0,\ldots,l), q_k\;(k=0,\ldots,m)$ polynomials) 
are different, and assume without loss of generality that $l>m$. 
Then, differentiate {\tt de2} $l-m$ times to get equations
\[
S_p:=\sum_{k=0}^{p} q_k^p\, g_n^{(k)}=0\quad\quad\quad(p=m,\ldots,l)
\;.
\]
Check if there are nontrivial rational functions 
$A_p\not\equiv 0\;(p=m,\ldots,l)$ 
such that a linear combination $\sum\limits_{p=m}^l A_p\,S_p$
is equivalent to the left hand side of {\tt de1}, i.\ e.\ is a rational
multiple of it.

If this is not the case, then $f_n$ and $g_n$ do not satisfy a common
simple differential equation, and therefore are not identical; 
return this, and quit.
Otherwise they satisfy a common simple differential equation;
continue with (e).
\item[(e)]
Let $l$ be the order of the common simple differential equation for $f_n$
and $g_n$. For $k=0,\ldots,l-1$ check if $f_n^{(k)}(0)=g_n^{(k)}(0)$. 
(Note that by the holonomic structure
the knowledge of the initial values (\ref{eq:holonomicIV}) is sufficient
to generate those.) These initial conditions may depend on $n$,
and are proved by application of a discrete version of the same
algorithm. If one of these equations is falsified,
then the identity $f_n\equiv g_n$ is disproved; return this, and quit. 
Otherwise, if all equations 
are verified, the identity $f_n\equiv g_n$ is proved.
\end{enumerate}
}
\end{algorithm}
\pr
By a well-known result about differential equations of the type considered,
the solution of an initial value problem 
\[
\sum_{k=0}^l p_k(x)\,f_n^{(k)}(x)=0
\;,
\quad\quad
f_n^{(k)}(0)=a_k\;(k=0,\ldots,l-1)
\]
is unique. To prove that $f_n$ and $g_n$ are identical, it therefore suffices 
to show that they satisfy a common differential equation, and the same 
initial values. This is done by our algorithm.
\eop
For the example expressions
\[
f_n(x):=
L_n^{(-1/2)}(x)
\]
and
\[
g_n(x):=
\frac{(-1)^n}{n!\,2^{2n}}\,H_{2n}\lk\sqrt x\rk
\]
we get the common differential equation
\[
2\, n\, f(x) + (1 - 2 x)\, f'(x) + 2\, x\, f''(x)=0
\;.
\]
Therefore to prove the identity
\[
L_n^{(-1/2)}(x)=
\frac{(-1)^n}{n!\,2^{2n}}\,H_{2n}\lk\sqrt x\rk
\;,
\]
(see e.\ g.\ \cite{AS}, (22.5.38)), it is enough to verify the
two initial equations $f_n(0)=g_n(0)$ and $f_n'(0)=g_n'(0)$. To establish 
the first of these conditions, with {\sc Mathematica}, e.\ g., we get

{\small
\begin{verbatim}
In[25]:= eq = Limit[LaguerreL[n,-1/2,x],x->0]==
         Limit[(-1)^n/(n!*2^(2*n))*HermiteH[2*n,Sqrt[x]],x->0]

                             1
         Pochhammer[1 + n, -(-)]        n
                             2      (-1)  Sqrt[Pi]
Out[25]= ----------------------- == ---------------
                Sqrt[Pi]                     1
                                    n! Gamma[- - n]
                                             2
\end{verbatim}
}\noindent
which is to be verified. In this situation, we establish the first order
recurrence equations for both sides

{\small
\begin{verbatim}
In[26]:= FindRecursion[Limit[LaguerreL[n,-1/2,x],x->0],n]

Out[26]= (-1 + 2 n) a[-1 + n] - 2 n a[n] == 0

In[27]:= FindRecursion[Limit[(-1)^n/(n!*4^n)*HermiteH[2*n,Sqrt[x]],x->0],n]

Out[27]= (-1 + 2 n) a[-1 + n] - 2 n a[n] == 0
\end{verbatim}
}\noindent
that coincide, so that it remains to prove the initial statement

{\small
\begin{verbatim}
In[28]:= eq /. n->0

Out[28]= True
\end{verbatim}
}\noindent
and we are done. Similarly one may prove the second initial value statement
$f_n'(0)=g_n'(0)$.

Applying the same method, 
(\ref{eq:difference differential equation}) can be proved by the calculations

{\small
\begin{verbatim}
In[29]:= SimpleDE[(n+1)*Bateman[n+1,x]-(n-1)*Bateman[n-1,x],x]

                         2          2    3                      2
Out[29]= (2 n - 2 x + 4 n  x - 4 n x  + x ) F[x] + (-2 n x + 2 x ) F'[x] +

                 2
>     (2 n - x) x  F''[x] == 0

In[30]:= SimpleDE[2*x*D[Bateman[n,x],x],x]

                         2          2    3                      2
Out[30]= (2 n - 2 x + 4 n  x - 4 n x  + x ) F[x] + (-2 n x + 2 x ) F'[x] +

                 2
>     (2 n - x) x  F''[x] == 0
\end{verbatim}
}\noindent
and using the initial values 
$F_n(0)=0$ and $F_n'(0)=-2$ (see \cite{KS1}, (11)).

Also, one can prove Clausen's formula
\be
\lk
_{2}F_{1}\left.\lk\begin{array}{c}
\multicolumn{1}{c}{\begin{array}{cc}a&b\end{array}}\\[1mm]
\multicolumn{1}{c}{a\!+\!b\!+\!1/2}\\
            \end{array}\rb x\rk
\rk^2
=
\;_{3}F_{2}\left.\lk\begin{array}{c}
\multicolumn{1}{c}{\begin{array}{ccc}2a&2b&a+b\end{array}}\\[1mm]
\multicolumn{1}{c}{\begin{array}{cc}a\!+\!b\!+\!1/2&2a\!+\!2b\end{array}}
            \end{array}\rb x\rk
\;,
\label{eq:Clausen's formula}
\ee
generating the common differential equation
\beao
&&  8\,a\,b\,\left( a + b \right) \,f(x) +
\\&&
    2\,( -a - 2\,{a^2} - b - 4\,a\,b - 2\,{b^2} + x + 3\,a\,x +
       2\,{a^2}\,x + 3\,b\,x + 8\,a\,b\,x + 2\,{b^2}\,x ) \,f'(x) +
\\&&
    3\,x\,\left( -1 - 2\,a - 2\,b + 2\,x + 2\,a\,x + 2\,b\,x \right) \, f''(x) +
\\&& 2\,\left( -1 + x \right) \,{x^2}\,f'''(x) = 0
\eeao
for both sides of (\ref{eq:Clausen's formula}),
or other hypergeometric identities like the Kummer transformation
\[
_1 F_1\left.\lk\begin{array}{c}
a\\b \end{array}\rb x\rk
=e^x\;
_1 F_1\left.\lk\begin{array}{c}
b-a\\b \end{array}\rb -x\rk
\]
or like 
\[
_{0}{F}_{1}\left.\lk\begin{array}{c}
a \!\end{array}\rb x\rk
\cdot
\;_{0}{F}_{1}\left.\lk\begin{array}{c}
b \!\end{array}\rb x\rk
=
\;_{2}{F}_{3}\left.\lk\begin{array}{c}
\multicolumn{1}{c}{\begin{array}{cc}
\frac{a+b}{2}&\frac{a+b-1}{2}\end{array}}\\[1mm]
\multicolumn{1}{c}{\begin{array}{ccc}
a&b&a+b-1
\end{array}}
\end{array}\rb 4\,x\rk
\]
and
\[
_{1}{F}_{1}\left.\lk\begin{array}{c}
a\\ b \end{array}\rb x\rk
\cdot
\;_{1}{F}_{1}\left.\lk\begin{array}{c}
a\\ b \end{array}\rb -x\rk
=
\;_{2}{F}_{3}\left.\lk\begin{array}{ccc}
\multicolumn{1}{c}{\begin{array}{cc}
a&b-a \end{array}}\\[1mm]
\multicolumn{1}{c}{\begin{array}{ccc}
b&\frac{b}{2}&\frac{b+1}{2} 
\end{array}}
\end{array}\rb \frac{x^2}{4} \rk
\]
corresponding to the Kummer differential equation
\[
a\,f(x) - (b - x)\,f'(x) - x\, f''(x) = 0
\;,
\]
and to
\beao
&&
  \left( 1 - a - b \right) \,\left( a + b \right) \,f(x) +
    \left( - a\,b + {a^2}\,b + a\,{b^2} - 2\,x - 4\,a\,x -
       4\,b\,x \right) \,f'(x) +
\\&&+
    \left( a + {a^2} + b + 3\,a\,b + {b^2} - 4\,x \right) \,x\,f''(x) +
    2\,\left( 1 + a + b \right) \,{x^2}\,f'''(x) + {x^3}\,f''''(x) = 0
\;,
\eeao
and
\beao
&&
4\,a\,\left( a - b \right) \,x\,f(x) +
    \left( b - 3\,{b^2} + 2\,{b^3} - {x^2} - 2\,b\,{x^2} \right) \,f'(x) 
\\&&+
    x\,\left( -b + 5\,{b^2} - {x^2} \right) \,f''(x) +
    \left( 1 + 4\,b \right) \,{x^2}\,f'''(x) + {x^3}\,f''''(x) = 0
\;,
\eeao
respectively.

Note that one can also reverse the order of the algorithm, 
i.\ e.\ first find common recurrence 
equations for $f_n$ and $g_n$ with respect to $n$, and then check the initial
conditions (depending on $x$) with the aid of differential equations.
This method should be compared with recent results of Zeilberger 
(\cite{Zei1}--\cite{Zei3}).

Moreover the given algorithm is easily extended to the case of several 
variables, if the family given forms an admissible family with respect to all 
of its variables, i.\ e.\ for each variable exists
\begin{itemize}
\item[-]
either a simple recurrence equation (corresponding to a ``discrete'' variable),
\item[-]
or a simple derivative rule (corresponding to a ``continuous'' variable),
depending on shifts with respect to one of the discrete variables.
\end{itemize}
Note, however, that (for the moment) the algorithm only works if $f$ and $g$ 
are ``expressions'', and no symbolic sums, derivatives of symbolic order,
etc.\ occur. In the next sections, we will, however, extend the
above algorithm to these situations.

\section{Algorithmic verification of Rodrigues type formulas}
\label{sec:Algorithmic verification of Rodrigues type formulas}

Here we present an algorithm to verify identities of the Rodrigues type
\[
g(n,x)=f^{(n)}(n,x)
\quad\quad\quad \;(f, g\;\mbox{functions}\;,\quad n\;\mbox{symbolic})
\;.
\]
This algorithm, however, does only work if the function $f$ is of the
hypergeometric type. On the other hand, for most Rodrigues type formulas 
in the literature, see e.\ g.\ \cite{AS}, this condition is valid.

The procedure is based on the following 
\begin{algorithm}
{\bf (Find differential equation for derivatives of symbolic order)}
\label{algo:Rodrigues}
{\rm
Let $f$ be of the hypergeometric type, i.\ e.\ there is a Laurent-Puiseux
type representation $f(n,x)=\sum_k a_k x^k$. Then there is a simple differential
equation for $g(n,x):=f^{(n)}(n,x)$ 
which can be obtained by the following algorithm:
\begin{enumerate}
\item[(a)]
{\tt de1:=SimpleDE(f,x)}: \\
Calculate the simple differential equation 
{\tt de1} of $f$, see Algorithm~\ref{algo:Find a simple DE}.
\item[(b)]
{\tt re1:=DEtoRE(de1,f,x,a,k)}: \\
Transfer the differential equation {\tt de1}
into the corresponding recurrence equation {\tt re1} for $a_k$, see 
\cite{Koe92}, \S 6.
\item[(c)]
If {\tt re1} is not of the hypergeometric type
(or is not equivalent to the hypergeometric type \cite{P}),
then quit.
\item[(d)]
{\tt re2:=SymbolicDerivativeRE(re1,a,k,n)}: \\
Otherwise set $c_k:=(k+1)_n\,a_{k+n}$. Bring {\tt re1} into the form
\[
a_{k+m}=R(k)\,a_k
\;,
\]
rational $R$, and calculate the hypergeometric type recurrence equation
{\tt re2}
\be
c_{k+m}=\frac{(k+n+1)_m}{(k+1)_m}\,R(k+n)\,c_k
\label{eq:REdiffN}
\ee
for $c_k$.
\item[(e)]
{\tt de2:=REtoDE(re2,a,k,G,x)}:\\
Transfer the recurrence equation {\tt re2} into the corresponding  
differential equation {\tt de2} for the $n$th derivative 
$g(n,x):=f^{(n)}(x)$ of $f$, see \cite{Koe92}, \S~11.
\end{enumerate}
}
\end{algorithm}
\pr
Parts (a), (b) and (e) of the algorithm are described precisely in \cite{Koe92}.
Now, assume, $g(n,x)=f^{(n)}(n,x)$, and that $f$ has the representation 
$f(n,x)=\sum_k a_k x^k$. Then we get
\[
\sum_{k}  c_k\,x^k=g(n,x)=
f^{(n)}(n,x)=\sum_{k} (k+1-n)_n\,a_k\,x^{k-n}
=\sum_{k} (k+1)_n\,a_{k+n}\,x^{k}
\;.
\]
Therefore we have $c_k=(k+1)_n\,a_{k+n}$, and we get the recurrence
equation
\beao
c_{k+m}&=&
(k+m+1)_n\,a_{k+n+m}=
(k+m+1)_n\,R(k+n)\,a_{k+n}
\\[1mm]&=&
\frac{(k+m+1)_n}{(k+1)_n}\,R(k+n)\,c_k
=
\frac{(k+m+n)!}{(k+m)!}\frac{k!}{(k+n)!}\,R(k+n)\,c_k
\\&=&
\frac{(k+m+n)!}{(k+n)}\frac{k!}{(k+m)!}\,R(k+n)\,c_k
=
\frac{(k+n+1)_m}{(k+1)_m}\,R(k+n)\,c_k
\;,
\eeao
and hence (\ref{eq:REdiffN}), for $c_k$. This finishes the proof.
\eop
As a first example we consider the identity
\[
\erfc_n(x)=\frac{(-1)^n\,e^{-x^2}}{2^n\,n!}\dedn{x}{n}\lk e^{x^2}\erfc x\rk
\]
(see e.\ g.\ \cite{AS}, (7.2.9)), or equivalently
\be
(-1)^n\,2^n\,n!\,e^{x^2}\,\erfc_n(x)=\dedn{x}{n}\lk e^{x^2}\erfc x\rk
\;.
\label{eq:erfcRodrigues}
\ee
Algorithm~\ref{algo:Rodrigues} yields step by step

{\small
\begin{verbatim}
In[31]:= de1=SimpleDE[E^(x^2)*Erfc[x],x]

Out[31]= -2 F[x] - 2 x F'[x] + F''[x] == 0

In[32]:= re1=DEtoRE[de1,F,x,a,k]

Out[32]= -2 (1 + k) a[k] + (1 + k) (2 + k) a[2 + k] == 0

In[33]:= re2=SymbolicDerivativeRE[re1,a,k,n]
                                           2
Out[33]= -2 (1 + k + n) a[k] + (2 + 3 k + k ) a[2 + k] == 0

In[34]:= de2=REtoDE[re2,a,k,G,x]

Out[34]= -2 (1 + n) G[x] - 2 x G'[x] + G''[x] == 0
\end{verbatim}
}\noindent
thus finally the differential equation
\[
-2\,(1 + n)\, g(x) - 2\, x\, g'(x) + g''(x)= 0
\]
for the function $\dedn{x}{n}\lk e^{x^2}\erfc x\rk$, 
which also can be obtained by the single statement

{\small
\begin{verbatim}
In[35]:= RodriguesDE[E^(x^2)*Erfc[x],x,n]

Out[35]= -2 (1 + n) F[x] - 2 x F'[x] + F''[x] == 0
\end{verbatim}
}\noindent
For the left hand term of (\ref{eq:erfcRodrigues}) we get

{\small
\begin{verbatim}
In[36]:= de3=SimpleDE[E^(x^2)*Erfc[n,x],x]

Out[36]= -2 (1 + n) F[x] - 2 x F'[x] + F''[x] == 0
\end{verbatim}
}\noindent
i.\ e.\ the same differential equation.

As next example we consider the Rodrigues type identity
(\ref{eq:Rodrigues Bateman})
for the Bateman functions, and rewrite it as
\be
\frac{n!\,e^{-x}}{x}\,F_n(x)
=
\frac{d^n}{dx^n}\lk e^{-2x}\,x^{n-1}\rk
\;.
\label{eq:Rodrigues Bateman2}
\ee
Our implementation yields

{\small
\begin{verbatim}
In[37]:= RodriguesDE[E^(-2x)*x^(n-1),x,n]

Out[37]= 2 (1 + n) F[x] + 2 (1 + x) F'[x] + x F''[x] == 0

In[38]:= SimpleDE[E^(-x)/x*Bateman[n,x],x]

Out[38]= 2 (1 + n) F[x] + 2 (1 + x) F'[x] + x F''[x] == 0
\end{verbatim}
}\noindent
Algorithm~\ref{algo:Rodrigues} shows the applicability of
Algorithm~\ref{algo:Verification of identities} if in the expressions
involved Rodrigues type expressions occur, as soon as we can handle the
initial values. Since in Algorithm~\ref{algo:Rodrigues} the function $f$
is assumed to be of hypergeometric type, this, however, can be done by
a series representation using 
Algorithm~\ref{algo:Find a Laurent-Puiseux representation} if $f$ moreover 
is analytic, and if the function $f$ of Algorithm~\ref{algo:Rodrigues}
does not depend on $n$: In this case
Algorithm~\ref{algo:Find a Laurent-Puiseux representation}
generates the generic coefficient $a_k$ of the series representation
$f(x)=\sum\limits_{k=0}^\infty a_k\,x^k$, and therefore we get the 
initial values by Taylor's theorem:
\[
\lk\dedn {x}{n} f\rk(0)=n!\,a_n\;.
\]
In our first example we conclude

{\small
\begin{verbatim}
In[39]:= PowerSeries[E^(x^2)*Erfc[x],x]

              2 k
             x
Out[39]= Sum[----, {k, 0, Infinity}] +
              k!

              k  1 + 2 k
          -2 4  x        k!
>    Sum[-------------------, {k, 0, Infinity}]
         Sqrt[Pi] (1 + 2 k)!
\end{verbatim}
}\noindent
so that the first initial condition for identity (\ref{eq:erfcRodrigues})
is given by the calculation (see \cite{AS} (7.2.7))
\[
\frac{(-1)^n\,n!}{\Gamma\lk \frac{n}{2}+1\rk}=
(-1)^n\,2^n\,n!\,\erfc_n(0)=\dedn{x}{n}\lk e^{x^2}\erfc x\rk(0)
=
\funkdeff{\frac{(2k)!}{k!}}{n=2k\;(k\in\N_0)}
{-{{2\,k!\,{4^k}}\over {{\sqrt{\pi }}}}
}{n\!=\!2k\!+\!1\;(k\in\N_0)}
\!\!\!,
\]
and the second one is established similarly. 

To identify the first initial values of our second example, we proceed as
follows: The left hand side of (\ref{eq:Rodrigues Bateman2}) yields
\be
\lim_{x\pf 0}\frac{n!\,e^{-x}}{x}\,F_n(x)
=
n!\,\lim_{x\pf 0}\frac{F_n(x)}{x}=n!\,F_n'(0)=-2\,n!
\label{eq:agreement}
\ee
(see \cite{KS1}, (11)), 
whereas from the identity
\[
\Big( x^n\Big)^{(k)}(0)=\funkdef{n!}{k=n}{0}
\]
and Leibniz's formula we derive for the right hand side
\beao
\lk e^{-2x}\,x^{n-1}\rk^{(n)}(0)
&=&
\lk\sum_{k=0}^n \ueber{n}{k} \lk x^{n-1}\rk^{(k)}
\lk e^{-2x}\rk^{(n-k)}\rk(0)
\\&=&
\ueber{n}{n-1}\,(n-1)!\,\lk e^{-2x}\rk'(0)=
-2\,n!
\;,
\eeao
in agreement with (\ref{eq:agreement}).

It is easily seen that we can always identify the initial values 
algorithmically by the method given if
$f(n,x)=w(x)\,X(x)^n$ with a polynomial $X$, i.\ e.\ is of the form
(\ref{eq:Rodriguestype}).

These results are summarized by
\begin{algorithm}
{\bf (Verification of identities)}
{\rm
With Algorithms~\ref{algo:Verification of identities} and 
\ref{algo:Rodrigues} identities involving Rodrigues type expressions can be
verified if only symbolic derivatives $f^{(n)}$ of hypergeometric type analytic 
expressions $f$ occur that have the form $f(n,x)=w(x)\,X(x)^n$ for some
polynomial $X$.
}
\end{algorithm}

\section{Algorithmic verification of formulas involving symbolic sums}
\label{sec:Algorithmic verification of formulas involving symbolic sums}

In this section we study, how identities involving symbolic sums can
be established. The results depend on the following algorithm
(compare \cite{SZ}, {\sc Maple} function {\tt cauchyproduct}):

\begin{algorithm}
{\bf (Find recurrence equation for symbolic sums)}
\label{algo:symbolic sums}
{\rm
Let $f_n(x)$ form an admissible family, and let $s_n(x)$ denote the 
symbolic sum $s_n(x):=\sum\limits_{k=0}^n f_k(x)$. Then the following
algorithm generates a recurrence equation for $s_n$:
\begin{enumerate}
\item[(a)]
{\tt re:=FindRecursion(f,k)}: \\
Calculate the simple recurrence equation
{\tt re} of $f_k$, see \cite{Koe92}, \S 11.
\item[(b)]
{\tt de1:=REtoDE(re1,f,k,F,z)}: \\
Transfer the recurrence equation {\tt re}
into the corresponding differential equation {\tt de1} valid for the
generating function $F(z):=\sum\limits_{k=0}^\infty f_k(x)\,z^k$,
see \cite{Koe92}, \S 11.
\item[(c)]
{\tt de2:=F(z)+(z-1)*F'(z)=0}: \\
Let {\tt de2} be the differential equation corresponding to the function
\[
G(z):=\sum\limits_{k=0}^\infty g_k z^k=\sum\limits_{k=0}^\infty z^k=
\frac{1}{1-z}
\;.
\]
\item[(d)]
{\tt de:=ProductDE(de1,de2,F,z)}: \\
Calculate the simple differential equation {\tt de} corresponding to the 
product $H(z):=F(z)\,G(z)$, see \cite{KS}, Theorem 3 (d). This 
differential equation has the order of {\tt de1}.
\item[(e)]
{\tt re:=DEtoRE(de,F,z,s,n)}:\\
Transfer the differential equation {\tt de} into the corresponding  
recurrence equation {\tt re} for the coefficient $s_n$ of $H(z)$,
see \cite{Koe92}, \S~6.
\end{enumerate}
}
\end{algorithm}
\pr
Parts (a), (b) and (e) of the algorithm are described precisely in \cite{Koe92}.
The rest follows from the Cauchy product representation
\[
H(z)=F(z)\,G(z)=\sum_{n=0}^\infty \lk\sum_{k=0}^n f_k\,g_{n-k}\rk\,z^n
=\sum_{n=0}^\infty \lk\sum_{k=0}^n f_k\rk\,z^n
\]
of the product function $F(z)\,G(z)$.
\eop
As an example we consider the sum $\sum\limits_{k=0}^n L_{k}^{(\alpha)}(x)$.
We get stepwise:

{\small
\begin{verbatim}
In[40]:= re=FindRecursion[LaguerreL[k,alpha,x],k]

Out[40]= (-1 + alpha + k) a[-2 + k] + (1 - alpha - 2 k + x) a[-1 + k] +

>     k a[k] == 0

In[41]:= de1=REtoDE[re,a,k,F,z]

                                                       2
Out[41]= (-1 - alpha + x + z + alpha z) F[z] + (-1 + z)  F'[z] == 0

In[42]:= de2=F[z]+(z-1)*F'[z]==0;

In[43]:= de=ProductDE[de1,de2,F,z]

                                                         2
Out[43]= (-2 - alpha + x + 2 z + alpha z) F[z] + (-1 + z)  F'[z] == 0

In[44]:= DEtoRE[de,F,z,s,n]

Out[44]= (2 + alpha + n) s[n] + (-4 - alpha - 2 n + x) s[1 + n] + 

>     (2 + n) s[2 + n] == 0
\end{verbatim}
}\noindent
or by a single statement

{\small
\begin{verbatim}
In[45]:= re=SymbolicSumRE[LaguerreL[k,alpha,x],k,n]

Out[45]= (2 + alpha + n) a[n] + (-4 - alpha - 2 n + x) a[1 + n] +

>     (2 + n) a[2 + n] == 0
\end{verbatim}
}\noindent
and substituting $n$ by $n-2$

{\small
\begin{verbatim}
In[46]:= Simplify[re /. n->n-2]

Out[46]= (alpha + n) a[-2 + n] + (-alpha - 2 n + x) a[-1 + n] + n a[n] == 0
\end{verbatim}
}\noindent
On the other hand, the calculation

{\small
\begin{verbatim}
In[47]:= FindRecursion[LaguerreL[n,alpha+1,x],n]

Out[47]= (alpha + n) a[-2 + n] + (-alpha - 2 n + x) a[-1 + n] + n a[n] == 0
\end{verbatim}
}\noindent
shows that the left and right hand sides of the identity
\be
\sum_{k=0}^n L_k^{(\alpha)}(x)=L_n^{(\alpha+1)}(x)
\label{eq:example identity}
\ee
(see e.\ g.\ \cite{Tri}, VI (1.16)) satisfy the same recurrence equation.

In our example identity two initial values remain to be considered
\[
L_0^{(\alpha)}(x)=L_0^{(\alpha+1)}(x)=1
\quad\quad\mbox{and}\quad\quad
L_0^{(\alpha)}(x)+L_1^{(\alpha)}(x)=L_1^{(\alpha+1)}(x)=2 + \alpha - x
\]
that trivially are established. 

Thus Algorithm~\ref{algo:symbolic sums} shows the applicability of
Algorithm~\ref{algo:Verification of identities} if in the expressions
involved symbolic sums occur. This is summarized by

\begin{algorithm}
{\bf (Verification of identities)}
{\rm
With Algorithms~\ref{algo:Verification of identities} and 
\ref{algo:symbolic sums} identities involving symbolic sums can be
verified.
\hfill$\Box$
}
\end{algorithm}
We like to mention that the function {\tt FindRecursion} is successful
for composite $f_n$ 
as long as recurrence equations exist and are applied recursively.
Here obviously no derivative rules are needed. 

We note further that as a byproduct this
algorithm in an obvious way can be generalized to sums 
$\sum\limits_{k=0}^n a_k\,b_{n-k}$ of the Cauchy product type. As an
example, the algorithm generates the recurrence equation
\be
2 \,(1 + 2 n)\, s_n - (1 + n)\, s_{n+1} = 0
\label{eq:sumbinomial}
\ee
for $s_n:=\sum\limits_{k=0}^n \ueber{n}{k}^2=
n!^2 \sum\limits_{k=0}^n \frac{1}{k!^2}\frac{1}{(n-k)!^2}$,
{\small
\begin{verbatim}
In[48]:= re=ConvolutionRESum[1/k!^2,1/k!^2,k,n]

                                   3
Out[48]= 2 (1 + 2 n) a[n] - (1 + n)  a[1 + n] == 0

In[49]:= ProductRE[re,FindRecursion[n!^2,n],a,n]

Out[49]= 2 (1 + 2 n) a[n] + (-1 - n) a[1 + n] == 0
\end{verbatim}
}\noindent
compare \cite{Zei1}--\cite{Zei3}.

Algorithm~\ref{algo:symbolic sums} may further be used to find a closed form
representation of a symbolic sum in case the resulting term is
hypergeometric:

\begin{algorithm}
{\bf (Closed forms of hypergeometric symbolic sums)}
{\rm
Let $s_n:=\sum\limits_{k=0}^n f_k$ be a hypergeometric term, i.\ e.\
$\frac{s_{n+1}}{s_n}$ be a rational function, then the following
procedure generates a closed form representation for $s_n$:
\begin{enumerate}
\item[(a)]
{\tt re:=SymbolicSumRE(f,k,n)}: \\
Calculate the simple recurrence equation {\tt re} of $s_n$ using 
Algorithm~\ref{algo:symbolic sums}.
\item[(b)]
If {\tt re} is of the hypergeometric type, then solve it by the hypergeometric
coefficient formula, else
apply Petkovsek's algorithm to find the hypergeometric solution $s_n$
of {\tt re}.
\hfill$\Box$
\end{enumerate}
}
\end{algorithm}
This result should be compared with the Gosper algorithm \cite{Gos}. Our
procedure is an alternative decision procedure for the same purpose.
Note that from the hypergeometricity of $s_n$ the hypergeometricity of $f_k$
follows \cite{Gos}, so that the first step of Algorithm~\ref{algo:symbolic sums}
leads to a simple first order recurrence equation.

Applying our algorithm to our example case 
$s_n=\sum\limits_{k=0}^n \ueber{n}{k}^2$, we get from
(\ref{eq:sumbinomial}), and the initial value $s_0=1$ the representation
\[
s_n=4^n\,\frac{\lk \ed 2\rk_n}{n!}=\frac{(2n)!}{n!^2}=\ueber{2n}{n}
\;.
\]
On the other hand, for $s_n=\sum\limits_{k=0}^n \ueber{n}{k}^3$, our procedure
gives
{\small
\begin{verbatim}
In[50]:= re=ConvolutionRESum[1/k!^3,1/k!^3,k,n]

                                          2
Out[50]= 8 a[n] + (1 + n) (16 + 21 n + 7 n ) a[1 + n] -

                     5
>     (1 + n) (2 + n)  a[2 + n] == 0

In[51]:= ProductRE[re,FindRecursion[n!^3,n],a,n]

                  2                        2                    2
Out[51]= 8 (1 + n)  a[n] + (16 + 21 n + 7 n ) a[1 + n] - (2 + n)  a[2 + n] ==

>    0
\end{verbatim}
}\noindent
and by Petkovsek's algorithm it turns out that $s_n$ is no hypergeometric term.

%
%
%
%
%

\end{document}